\documentclass[11pt]{article}
\usepackage{graphics}
\usepackage{amsmath}
\usepackage{amssymb}
\usepackage{amsthm}
\usepackage{mathtools}
\usepackage{glossaries}
\usepackage{longtable}
\usepackage{lscape}
\usepackage{rotating}
\usepackage{bm}
\usepackage{hyperref}
\usepackage{color}
\usepackage{comment}
\usepackage{enumitem}
\usepackage{caption}
\newcommand{\myemph}[1]{\textbf{#1}}
\newcommand{\klemme}{\vspace{-1.5mm}}

\newcommand{\xa}{\overline{a}}        
\newcommand{\xb}{\overline{b}}        
\newcommand{\xz}{\lambda}             
\newcommand{\xk}{\kappa}              
\newcommand{\xM}{\mbox{$\cal{T}$}}    
\newcommand{\xNT}{\mbox{$\cal{N}$}}   
\newcommand{\xeuro}{\texteuro}        

\title{Decision support for the Technician Routing and Scheduling Problem}
\author{Mette Gamst, David Pisinger}
\date{\small DTU Management, Technical University of Denmark, Kgs. Lyngby, Denmark}

\begin{document}
\setlist{nosep}

\maketitle

\begin{small}
\noindent\textbf{Abstract}
The technician routing and scheduling problem (TRSP) consists of technicians serving tasks subject to qualifications, time constraints and routing costs. In the literature, the TRSP is solved either to provide actual technician plans or for performing what-if analyses on different TRSP scenarios. We present a method for building optimal TRSP scenarios, e.g., how many technicians to employ,  which technician qualifications to upgrade, etc. The scenarios are built such that the combined TRSP costs (OPEX) and investment costs (CAPEX) are minimized.  
Using a holistic approach we can generate scenarios that would not have been found by studying the investments individually.
The proposed method consists of a matheuristic based on column generation. To reduce computational time, the routing costs of a technician are approximated.
The proposed method is evaluated on data from the literature and on real-life data from a telecommunication company. The evaluation shows that the proposed method successfully suggests attractive scenarios. The method especially excels in ensuring that more tasks are serviced but also reduces travel time with around 16\% in the real-life instance. 
We believe that the proposed method could constitute an important strategic tool in field service companies and we propose future research directions to further its applicability.
\end{small}



\section{Introduction}
The \emph{Technician Routing and Scheduling problem} (TRSP) consists of assigning tasks to technicians subject to skill sets, time windows, travel times and routing costs. The problem occurs in telecommunications, where technicians perform installation or fault correction tasks. It is closely related to more general workforce scheduling problems such as general maintenance scheduling, home care scheduling, etc. \cite{Bruecker, Pereira}. The TRSP and workforce scheduling problems have been extensively studied in the research literature, however, mainly with an operational focus on generating actual work schedules \cite{Salazar}. 

This paper considers the strategic decisions of determining the right set of technicians, technician qualifications, task outsourcing or digitization, and working hours such that the operational costs can be minimized. The strategic decisions all come with an investment cost. Traditionally, strategic questions are analyzed via what-if analyses. Given a base case and a \emph{scenario}, e.g., with a distinct set of technicians, the difference in the operational costs indicates whether the scenario is an attractive investment \cite{Kottemann}. We propose a method which minimizes the combined investment costs and operational costs (CAPEX + OPEX) to generate attractive scenarios. To the best of our knowledge, this has not been proposed previously for the TRSP. In a deterministic and optimal setting, this replaces the traditional what-if approach. In practice, however, the approach can be viewed as a scenario builder, where the scenarios afterwards can be evaluated in detail using what-if analysis.
By considering all investments in a holistic model, we are able to study scenarios that would
not have been found by traditional what-if tools that analyze investments one by one.

To reach computational tractability, the proposed method approximates the TRSP by assigning tasks to technicians instead of optimizing each technician path. We show that this is a fair approximation by comparing it to the TRSP solved by an \emph{Adaptive Large Neighborhood Search} metaheuristic. The assignment problem is solved through column generation, where the master problem handles the investments and the subproblem assigns tasks to each technician. The method is computationally evaluated on data from the literature and on real-life TRSP data. 

The main contributions of this paper are:
\begin{itemize}
    \item A general method to optimize investments in the TRSP using a holistic approach that combines CAPEX and OPEX 
    \item Evaluation of investments in the technician fleet, overtime, skill sets, and digitization of tasks on benchmark data and real-life TRSP data
    \item An Adaptive Large Neighborhood Search for the TRSP
    \item An approximation of the TRSP based on assigning tasks to technicians
    \item Publication of a real-life large TRSP dataset
\end{itemize}

\noindent
The paper is organized as follows: The TRSP is formally defined in Section~\ref{problem_definition}, and related literature is reviewed in Section~\ref{literature}. Next, Section~\ref{tacttrsp} presents the proposed solution method, approximating TRSP with a task assignment approach. The assignment approach is extended with investments in Section~\ref{tactinv}. Section~\ref{data} presents test data, Section~\ref{computations} computational experiments, and Section~\ref{future} future work. Finally, Section~\ref{conclusion} concludes the paper.

\section{Problem definition}
\label{problem_definition}
Let $\mathcal{T}$ be a set of \emph{technicians}, $V$ the set of \emph{tasks}, and $S$ the set of all \emph{skills}. Each technician $t \in \mathcal{T}$ has a \emph{start} and \emph{endpoint} (depot) $d_t$, \emph{skills} set $S(t) \subseteq S$ and a set of \emph{work shifts} $D(t)$ consisting of \emph{time windows} $[\xa_{t}^{d}, \xb_{t}^{d}]$. 
Each task $i \in V$ consists of the \emph{required skills} set $S(i)$, a \emph{duration} $f_i$, 
a \emph{time window} $[a_{i}, b_{i}]$, 
and a \emph{penalty} for not being serviced $c_i$.

We generate a \emph{daily technician} for each workday for each technician in the set $\mathcal{T}$. Let the set of daily technicians be denoted $T$. Each daily technician $t\in T$ has a \emph{start} and \emph{endpoint} (depot) $d_t$, \emph{skills} set $S(t)$ and \emph{time window} $[\xa_{t}, \xb_{t}]$. The set of tasks, that can be served by a daily technician $t \in T$ subject to skills and time windows, is denoted $V(t)$.
The \emph{travel time} between $i$ and $j$ is denoted $c_{ij}^t$ for technician $t$. Travel times are assumed to satisfy the triangle inequality.
We introduce binary variables $x_{ij}^t$ to denote if technician $t$ travels from task $i$ to $j$, binary variables $y_i$ to denote if
task $i$ is unassigned, and linear variables $z_i^t$ to denote the arrival time at task $i$.
 
The TRSP consists of minimizing the time and costs associated with servicing the tasks subject to skill constraints and time windows. 
It can be now be modelled as:

\begin{small}
\begin{align}
    \min \quad & \sum_{t \in T}\sum_{i \in V(t)}\sum_{j \in V(t)} c_{ij}^t x_{ij}^t + \sum_{i \in V(t)} c_i y_i \label{trsp_obj} \\
    \text{s.t.} \quad & \sum_{\begin{subarray}{c}t\in T:\\i \in V(t)\end{subarray}} \sum_{j \in V(t)} x_{ij}^t + y_i \geq 1 & \forall i \in V \label{trsp_c1}\\
    & \sum_{i\in V(t)} x_{d_ti}^t \leq 1 & \forall t\in T\label{trsp_c2}\\
    & \sum_{j\in V(t)} x_{ij}^t - \sum_{j\in V(t)} x_{ji}^t = 0 & \forall t \in T, \forall i\in V(t) \label{trsp_c3}\\
    & \sum_{i\in V(t)} x_{id_t}^t \leq 1 & \forall t\in T\label{trsp_c4}\\
    & z_i^t + f_i + c_{ij}^t \leq z_j^t + M (1-x_{ij}^t) & \forall t \in T, \forall i, j \in V(t) \label{trsp_c5}\\
    & a_i \leq z_i^t \leq b_i-f_i & \forall t \in T, \forall i \in V(t) \label{trsp_c6}\\
        & \xa_t \leq z_i^t \leq \xb_t-f_i-c_{id_t}^t & \forall t \in T, \forall i \in V(t) \label{trsp_c7}\\
& x_{ij}^t \in \{0, 1\} & \forall t\in T, \forall i, j \in V(t) \label{trsp_b1}\\
& y_i \geq 0 & \forall i \in V(t) \label{trsp_b1a}\\
& z_i^t \geq 0 & \forall t \in T, \forall i \in V(t) \label{trsp_b2}
\end{align}
\end{small}
The objective function \eqref{trsp_obj} minimizes the total travel time and the total penalty for unserved tasks. The penalty is given as a measure of time to make the two terms in the objective function comparable. The first constraints \eqref{trsp_c1} ensure that a task is either served or marked as unserved. The next three constraint sets \eqref{trsp_c2}-\eqref{trsp_c4} are flow conservation constraints saying that a daily technician can leave home at most once, must leave any entered task and return to home at most once.  The final three constraint sets \eqref{trsp_c5}-\eqref{trsp_c7} ensure time windows; constraints \eqref{trsp_c5} keep track of start times at each visited task for each daily technician. Constraints \eqref{trsp_c6}-\eqref{trsp_c7} ensure that task start times are within both the time windows of the task and of the technician. Finally, bounds \eqref{trsp_b1}-\eqref{trsp_b2} ensure that variables take on feasible values. Note that the bound on variables $y_i$ can be relaxed without any loss of information.

It is noted that the formulation bears strong resemblance with the \emph{Multi-Depot Vehicle Routing Problem with Time Windows} \cite{laporte1984} and only differs in the definition of the set $V(t), t\in T$ and that tasks may be unserved. 


\section{Literature review}
\label{literature}
The TRSP was introduced by Tsang and Voudouris \cite{tsang} with an application in British Telecom. They solved the TRSP using a guided local search. The TRSP has since been extensively studied and mainly solved using metaheuristics \cite{mathlouthi_tabu, pillac_matheur}. Few contributions attempt to solve the problem to optimality \cite{mathlouthi_mip,mathlouthi_bap} for small instances with at most 45 customers. In comparison most of the heuristic approaches consider instances with up to 200 customers. We refer to Chen et al. \cite{chen} and Mathlouthi et al. \cite{mathlouthi_tabu} for recent literature overviews. Most literature contributions consider real-life applications including specific extra constraints. Pillac et al. \cite{pillac_dynamic} consider the dynamic TRSP where new customers may arrive during the planned period. They propose an Adaptive Large Neighborhood Search algorithm and solve instances with up to 200 customers and 25 technicians. Mathlouthi et al. \cite{mathlouthi_tabu} consider a similar problem and apply Tabu Search with integrated adaptive memory. They generate a set of elite solutions which are stored in memory and used to construct new starting solutions for the tabu search.  The approach has also successfully been applied to the classic vehicle routing problem \cite{gendreau, vidal}. In another variant of the TRSP, service times may differ depending on how experienced a technician is, and that technicians may become more and more experienced and thus faster. Chen et al. \cite{chen} take this into account and solve the problem heuristically. To the best of our knowledge, the TRSP research literature only considers the problem of generating work schedules and/or technician teams, while the more strategic decisions such as the number and location of technicians or their skills set, have yet to be explored.

In recent years, the TRSP and related scheduling problems such as home care scheduling, have been considered in more overall terms as workforce scheduling problems. Recent surveys consider these closely related problems to promote learning across the problem types and to emphasize similarities in solution approaches and problem instance types \cite{Bruecker, Pereira}. 

The \emph{Vehicle Routing Problem} (VRP) is related to the TRSP and is an extensively studied problem. Relevant VRP variants are the multi-depot VRP \cite{md_hf_vrp}, the heterogeneous VRP \cite{md_hf_vrp} where the vehicles differ, the site dependent VRP where each customer can only be served by a subset of the vehicles \cite{golden_site} and the Vehicle Fleet Mix Problem where the usage of each vehicle is associated with a cost \cite{vmp}.  A survey on the VRP and its variants can be found in \cite{vrp_survey}. The Vehicle Fleet Mix Problem takes on a more strategic view by combining investment costs of using a vehicle with the operational costs of serving customers \cite{golden_fmp}. Many mathematical formulations and tailored heuristics have been proposed for the problem, see e.g. \cite{vmp,li_vmp,salhi_vmp}, but they are not easily extended to other strategic investment decisions. Other interesting work on VRP related problems includes Gaudioso and Paletta \cite{gaudioso} who suggest an alternative heuristic for the tactical problem of minimizing fleet size, rather than the operational problem of reducing travel costs. 
Moreover, the research literature focuses on making actual vehicle plans.
Another interesting research contribution is the work of Mourgaya and Vanderbeck \cite{vanderbeck} who solve a tactical version of the Periodic VRP in which visit schedules and customer assignments to vehicles are solved simultaneously. The sequencing of customers within vehicle paths is determined in an operational problem. The authors consider two objectives: a “workload balancing” objective that ensures an equal distribution of customers among vehicles and a “regionalization” objective that clusters customers geographically as a proxy for tour length. Focusing solely on the tactical problem, facilitates the solution of larger problem instances. 

The \emph{Facility Location Problem} combines the cost of opening facilities with the cost of serving customers \cite{fac_book}. This bears resemblance to combining investments in e.g., technicians with the routing costs in the TRSP. 
Real-life facility location applications include opening of factories, hospitals and placement of electric vehicle chargers \cite{fac_casestudy}. Numerous solution methods have been proposed for solving the Facility Location Problem, including MIP solvers, Lagrangian relaxation, Benders decomposition and a variety of heuristics, see e.g., the survey on Facility Location Problems in health care \cite{fac_health}. The Facility Location Problem differs from the investment problem in TRSP in the lack of routing and that only one type of investment is possible; the opening of a facility.

The idea of minimizing the combined CAPEX and OPEX is also present in other research areas. In energy analyses, the Capacity Expansion Problem (CEP) consists of investing in energy producing units subject to meeting energy demand as cheaply as possible \cite{cep_ori}. For example, how much renewable energy to invest in \cite{cep_res}. The CEP is NP-hard and often considers at least a yearly time horizon in hourly resolution in order to determine the best energy mix \cite{buchholz}. Attempts to achieve computational tractability includes Dantzig-Wolfe decomposition \cite{cep_dw}, Benders decomposition \cite{cep_ben}, and heuristics consisting of aggregating the time domain  \cite{buchholz,cep_aggr}.
A popular approach for aggregating the time domain, consists of clustering days to collapse the time horizon and thus to reduce the problem size into tractable MIPs. 

Considering both CAPEX and OPEX is also relevant in design problems such as \emph{Wind Farm Layout Problems}. Here the installation costs are considered together with wind production performance which depends on the wind farm layout and resulting wake effects \cite{fischetti}.

Combining CAPEX and OPEX in optimization is thus far from new. However, it does not seem to have gained much attention in the context of routing or the TRSP. This paper seeks to remedy this research gap.


\section{Task assignment for approximating the Technical Routing and Scheduling Problem}
\label{tacttrsp}
The TRSP is NP-hard and is typically solved in the literature using meta-heuristics such as local neighborhood search algorithms where the neighborhood consists of some method to swap tasks between technicians or permuting technician paths.  It is not straightforward to include investments in this approach for several reasons: 
\begin{itemize}
\item An investment is often only attractive if it results in many tasks being assigned to a technician. The neighborhood search must thus be extended to considering groups of tasks instead of individual tasks, which again may result in very large neighborhoods and poor convergence. 
\item An investment is often only attractive if it is utilized over a larger time period, e.g. months or even years. The investment cost can be amortized into daily costs, but the investment business case is very fragile if it is built on the TRSP for a single day. The TRSP can be solved for larger time periods using the mentioned local search heuristics, but then the neighborhood search must be extended to considering the path for a technician for each day in the entire time horizon, to investigate the full potential of the investment. This is not straightforward.
\end{itemize}
Generally, heuristics are  greedy and shortsighted in their nature and investments are only worthwhile considering large time horizons. Hence, we instead define a mathematical model for the investment problem.  To reduce problem complexity, we propose an assignment of tasks to technicians, similar to the tactical assignment approach of Mourgaya and Vanderbeck \cite{vanderbeck}. This is a simplification of solving the full TRSP, because travel times are estimated instead of optimized. 

The remainder of this section first studies the task assignment approach and afterwards an Adaptive Large Neighborhood Search algorithm for the full TRSP is presented. We need the latter to be able to evaluate the performance of the assignment approach. 

\subsection{Task assignment}
The task assignment problem consists of assigning tasks to technicians subject to skills and the time windows of the task and technician. The task assignment problem does not optimize each technician path and thus does not check for all overlapping service or tasks, nor is the exact travel time computed. The travel time is estimated by considering the distance between a technician's home depot and the location of assigned tasks.

The task assignment problem is defined using the same set and parameter notation as in Section~\ref{problem_definition}. Recall that a \emph{daily technician} is generated for each working day for each technician. In the following, we denote the original technician as the \emph{master technician}. 
The capacity $Q_t$ (i.e. time available) of technician $t$ is $\xb_{t}-\xa_{t}$.
Furthermore, we introduce a constant $k \geq 0$  to scale the estimated travel time such that it better reflects reality.

Variable $x_i^t \in \{0, 1\}$ indicates if task $i \in V$ is assigned to daily technician $t \in T$. Variable $y_i \in \{0, 1\}$ decides if task $i \in V$ is not assigned to any daily technician. Finally, variable $\xz_t \geq 0$ is an estimate on how far daily technician $t$ travels provided the current assignment. The task assignment problem is:
\begin{small}
\begin{align}
     \min \quad&  \sum_{t \in T} \xz_t + \sum_{i \in V} c_i y_i \label{full_obj}\\
     \text{s.t.} \quad& \xz_t \geq \xk \cdot c_{ij} (x_i^t + x_j^t - 1) & \forall i, j \in V(t), t\in T \label{full_c1} \\
      & \xz_t \geq \xk \cdot c_{d_ti} x_i^t & \forall i \in V(t), t\in T \label{full_c2}\\
      & \sum_{i \in V(t)}f_i x_i^t \leq Q_t -   \sum_{i \in V} c_{d_ti}x_i^t & \forall t \in T \label{full_c3}\\
      & \sum_{\begin{subarray}{c}t \in T:\\ i \in V(t)\end{subarray}} x_i^t + y_i \geq 1 & \forall i \in V \label{full_c4}\\
      & x_i^t + x_j^t \leq 1 & \forall i,j\in V(t): overlap(i,j), \forall t \in T \label{full_c5}\\
      & x_i^t \in \{0, 1\} & \forall i \in V(t), \forall t\in T \label{full_b1}\\
      & y_i \geq 0 & \forall i \in V \label{full_b2}\\
      &\xz_t \geq 0 & \forall t\in T \label{full_b3}
\end{align}
\end{small}
The objective function \eqref{full_obj} minimizes the total estimated travel time and the penalty for unassigned tasks. The first two sets of constraints, \eqref{full_c1} and \eqref{full_c2}, set the estimated total travel time for a daily technician to be the longest travel time between any two assigned tasks or from the depot to an assigned task. The travel time is multiplied with a scaling factor $\xk$ to account for additional driving to intermediate tasks. Constraints \eqref{full_c3} ensure that the daily technician capacity is satisfied, where we estimate time spent on traveling to be the travel time from the depot to each assigned task. Constraints \eqref{full_c4} say that each task must either be assigned or marked as unassigned. Because of the non-negative costs in the objective function, a task will never both be assigned and marked as unassigned. The last set of constraints \eqref{full_c5} ensure that a technician is not assigned two tasks, which cannot be serviced by the same technician because of overlapping time windows. Bounds \eqref{full_b1}-\eqref{full_b3} ensure feasible variables values.

\subsection{An Adaptive Large Neighborhood Search for the TRSP}
The task assignment approach is a simplification of the full TRSP. To evaluate its quality, it is compared to the full TRSP. The TRSP is NP-hard and difficult to solve. As seen in Section~\ref{literature}, the literature suggest a wide variety of meta-heuristics. We apply an Adaptive Large Neighborhood Search algorithm (ALNS) to solve the TRSP as formulated in mathematical model \eqref{trsp_obj} - \eqref{trsp_b2}. ALNS is one of the best performing meta-heuristics on Vehicle Routing Problems \cite{ALNS}. The ALNS iteratively destroys and repairs the current solution. After each destroy and repair operation, the resulting solution is evaluated and if certain requirements are met, it becomes the new current solution. Our implementation of ALNS follows the template outlined by Ropke and Pisinger \cite{ALNS}. We evaluate a solution by using the record-to-record method suggested in \cite{santini}. The destroy methods are:
\begin{itemize}
    \item \myemph{Random destroy}: random tasks are removed from daily technician paths and marked as unassigned
    \item \myemph{Route destroy}: a random daily technician path is selected and all tasks unassigned
\end{itemize}
The repair methods are:
\begin{itemize}
    \item \myemph{Greedy repair}: the cost of inserting each unrouted task in any path is calculated and the cheapest task is inserted first
    \item \myemph{Regret repair}: the cost of the best insert and the second best insert is calculated for each unrouted task. The task with highest different between the best and second best insert is inserted first \cite{ALNS}
\end{itemize}
When a task has been inserted or removed from a path, the arrival times in the path are updated to limit waiting times.

As a final optimization, the ALNS finishes with solving a set-cover-like formulation on the generated paths. Let $P$ be the set of generated paths, each path with cost $c_p, p\in P$ and let decision variable $x_p \in \{0, 1\}$ indicate if path $p \in P$ is selected or not:
\begin{small}
\begin{align}
     \min \quad &  \sum_{p \in P}\sum_{t \in T} c_p x_p + \sum_{i \in V} c_i y_i  \label{alns_obj}\\
     \text{s.t.} \quad & \sum_{p \in P} x_p^t \leq 1 & \forall t\in T \label{alns_c1}\\
     & \sum_{p \in P} \delta_i^p x_p^t + y_i  \geq 1 & \forall i \in V \label{alns_c2}\\
      & x_p^t \in \{0, 1\} & \forall p \in P, \forall t\in T \label{alns_b1}\\
        & y_i \geq 0 & \forall i \in V \label{alns_b2}
\end{align}
\end{small}
The objective \eqref{alns_obj} minimizes the total costs. Constraints \eqref{alns_c1} ensure that at most one path is used per technician, constraints \eqref{alns_c2} say that a task is either served or marked as unserved, and bounds \eqref{alns_b1}-\eqref{alns_b2} ensure that variables take on feasible values.

The number of paths $P$ may be very large. We only consider the 25000 most promising paths, where the potential of a path is defined by the value of the best ALNS solution it occurred in.

\section{Investment in the task assignment approach}
\label{tactinv}
To optimize investments, we minimize the combined OPEX and CAPEX costs. The task assignment problem is extended with investments by including the investment costs in the objective function and by adding investment constraints. First, we need to introduce some additional notation and define the exact investments to consider. Then follows the mathematical model for the problem.

Recall that the set $T$ contains daily technicians. An investment may, however, affect the \emph{master technician} and thus several daily technicians. Let $\xM$ be the set of master technicians and let parameter $\delta_m^t \in \{0, 1\}$ denote if $m\in \xM$ is the master technician for $t \in T$.

Let $A = \{ot,dig,s,nt\}$ be the set of possible investments, and variables $u_a^t \in \{0, 1\}, u_a^m \in \{0, 1\}, u_a^i \in \{0, 1\}$ indicate if investment $a \in A$ is performed, either for daily technician $t\in T$, master technician $m\in \xM$ or task $i \in V$. Each investment has an associated cost, $c_a^t$, $c_a^m$ resp. $c_a^i$. 

The considered investments are:
\begin{itemize}
    \item \myemph{Daily technician overtime}, i.e., if a daily technician should work overtime. Investment variable is $u_{ot}^t\in \{0, 1\}, t\in T$ and the available overtime minutes are given by parameter $w_{ot}$.
    \item \myemph{Digitization of tasks}. Instead of assigning a task, an investment can be made in the underlying equipment to enable remote fault handling. The investment variable is $u_{dig}^i\in \{0, 1\}, i\in V$ and parameter $\delta_{dig}^i \in \{0, 1\}$ indicate whether or not task $i \in V$ can be digitized.
    \item \myemph{Skill upgrade}, i.e., if the master technician should receive skill training. Investment variable is $u_{s}^m\in \{0, 1\}, m\in \xM$. The investment affects the skill ability of all corresponding daily technicians, i.e., all $t\in T: \delta_m^t = 1$.
    \item \myemph{New technicians}, i.e., if a new \emph{master technician} should be employed. The skills set, working days and location of the new master technician are given as input. Let $\xNT$ denote the corresponding set of new daily technicians and let the investment variable be denoted $u_{nt}^m \in \{0, 1\}$.
\end{itemize}
These are the investments, we consider in this paper. The list could be extended to any other investments, which would be handled in a similar manner. The mathematical model now becomes:
\begin{small}
\begin{align}
     \min \quad &  \sum_{t \in T} z_t + \sum_{i \in V} c_i y_i +
 \sum_{a \in A} \left(\sum_{t \in T} c_a^t u_a^t + \sum_{m \in \xM} c_a^m u_a^m + \sum_{i\in V} c_a^i u_a^i\right) \rule{-5cm}{0mm}\label{fullinv_obj}\\
     \text{s.t.} \quad & z_t \geq \xk \cdot c_{ij} (x_i^t + x_j^t - 1) & \forall i, j \in V, t\in T \label{fullinv_c1} \\
      & z_t \geq \xk \cdot c_{d_ti} x_i^t & \forall i \in V, t\in T \label{fullinv_c2}\\
    & \sum_{i \in V}f_i x_i^t \leq Q_t + w_{ot} \cdot u_{ot}^t -   \sum_{i \in V} c_{d_ti}x_i^t \rule{-5cm}{0mm}& \forall t \in T \label{fullinv_c3}\\
          & \sum_{t \in T} x_i^t + y_i + \delta_{dig}^i u_{dig}^i \geq 1 & \forall i \in V \label{fullinv_c4}\\
     & x_i^t + x_j^t \leq 1 & \forall t \in T, \forall i,j\in V: overlap(i,j) \label{fullinv_c4b}\\
      & x_i^t \leq  \sum_{m\in \xM} \delta_m^t u_s^m & \forall i \in V, \forall t\in T, \forall s \in S(i): s \not\in S(t) \label{fullinv_c5}\\
      & x_i^t \leq \sum_{m\in \xM} \delta_m^t u_{nt}^m & \forall i\in V, \forall t\in \xNT \label{fullinv_c6}\\
      & x_i^t = 0 & \forall i \in V, \forall t\in T: tw(i, t) = \emptyset \label{fullinv_c7}\\
      & x_i^t \in \{0, 1\} & \forall i \in V, \forall t\in T \label{fullinv_b1}\\
      & y_i \geq 0 & \forall i \in V \label{fullinv_b2}\\
      &z_t \geq 0 & \forall t\in T \label{fullinv_b3}\\
      & u_a^t \in \{0, 1\} & \forall a\in A, \forall t\in T\label{fullinv_b4}\\
      & u_a^m \in \{0, 1\} & \forall a\in A, \forall m\in \xM\label{fullinv_b5}\\
      & u_a^i \in \{0, 1\} & \forall a\in A, \forall i\in V\label{fullinv_b6}
\end{align}
\end{small}
The objective \eqref{fullinv_obj} minimizes the total estimated travel time, the penalties for unassigned tasks and the investment costs. We convert investment costs to equivalent time measurements to make the terms in the objective function comparable. 

Constraints \eqref{fullinv_c1}-\eqref{fullinv_c4b} correspond to \eqref{full_c1}-\eqref{full_c5}. The only extension is that constraints \eqref{fullinv_c3} satisfy time capacity for each daily technician subject to overtime investment. Furthermore, \eqref{fullinv_c4} ensure that a task is either assigned, marked as unassigned or digitized if possible. Because of the non-negative costs in the objective function, at most one of the three variable sets in the constraints will be set for each task.


The next two sets of constraints are new. Constraints \eqref{fullinv_c5} concern skill investment and state that if a required task skill, $s \in S(i), i\in V$ is not readily available by a daily technician $t\in T$, then the task can only be assigned to the daily technician if an investment is made in training the corresponding master technician $m\in \xM: \delta_m^t=1$. Constraints \eqref{fullinv_c6} concern investment in new technicians and state that a task can only be assigned to a new daily technician, if the corresponding investment is made in the corresponding master technician $m\in \xM: \delta_m^t=1$. The final constraints \eqref{fullinv_c7} enforce that a task cannot be assigned to a daily technician if their time windows do not overlap. Here, $tw(i, t)$ is the set of feasible start times that daily technician $t \in T$ can perform task $i \in V$ subject to task duration and time windows, including overtime. Bounds \eqref{fullinv_b1}-\eqref{fullinv_b6} ensure that variables are in the correct domain.

\subsection{Column generation}
The investment problem is decomposed into a master and subproblem. The master problem handles investment decisions and the subproblem assigns tasks to a daily technician. Let variable $x_p^t \in \{0, 1\}$ represent an assignment $p \in P$ for daily technician $t \in T$ with estimated travel time $c_p^t$. 
Furthermore, binary parameter $\delta_{ot}^p$ indicates if assignment $p$ uses overtime, $\delta_s^p$ if assignment $p$ requires skill $s$,
 and $\delta_i^p$ if assignment $p$ assigns task $i$. 
The master problem now becomes:

\begin{small}
\begin{align}
     \min \quad &  \sum_{p \in P}\sum_{t \in T} c_p^t x_p^t + \sum_{i \in V} c_i y_i  + \sum_{a \in A} \left(\sum_{t \in T} c_a^t u_a^t + \sum_{m \in \xM} c_a^m u_a^m + \sum_{i\in V} c_a^i u_a^i\right) \rule{-5cm}{0mm}\label{master_obj}\\
     \text{s.t.} \quad & \sum_{p \in P} x_p^t \leq 1 & \forall t\in T \label{master_c1}\\
     & \sum_{p \in P} x_p^t \leq \sum_{m \in \xM} \delta_m^t u_{nl}^m & \forall t\in \xNT \label{master_c2}\\
      & \sum_{p \in P} \delta_{ot}^p x_p^t \leq u_{ot}^t & \forall t \in T \label{master_c3}\\
      & \sum_{p \in P} \delta_{s}^p x_p^t \leq \sum_{m\in \xM} \delta_m^t u_{s}^m & \forall t \in T, \forall s \in S: s\not\in S(t) \label{master_c4}\\
      & \sum_{p \in P} \delta_i^p x_p^t + y_i + \delta_{dig}^i u_{dig}^i \geq 1 & \forall i \in V \label{master_c5}\\
      & x_p^t \in \{0, 1\} & \forall p \in P, \forall t\in T \label{master_b1}\\
      & y_i \geq 0 & \forall i \in V \label{master_b1a}\\
      & u_a^t \in \{0, 1\} & \forall a\in A, \forall t\in T\label{master_b2}\\
      & u_a^m \in \{0, 1\} & \forall a\in A, \forall m\in \xM\label{master_b3}\\
      & u_a^i \in \{0, 1\} & \forall a\in A, \forall i\in V\label{master_b4}
\end{align}
\end{small}
\noindent
The objective \eqref{master_obj} minimizes the total costs. The first set of constraints \eqref{master_c1} state that a daily technician can have at most one assignment. The second four constraint sets enforce investment costs. Constraints \eqref{master_c2} concern new technicians, \eqref{master_c3} overtime for daily technicians, \eqref{master_c4} skill upgrade and \eqref{master_c5} digitization of tasks together with ensuring that each task is either assigned, marked as unassigned or digitized. 
Bounds \eqref{master_b1}-\eqref{master_b4} ensure that variables take on feasible values.

%

The master problem is LP-relaxed and the reduced cost derived. The dual variables are $\pi_t^{\alpha} \leq 0, t\in T$ for constraints \eqref{master_c1}, $\pi_t^{\beta} \leq 0, t\in \xNT$ for constraints \eqref{master_c2}, $\pi_t^{\gamma} \leq 0, t\in T$ for constraints \eqref{master_c3}, $\pi_{st}^{\mu} \leq 0, t\in T, s\in S: s\not\in S(t)$ for constraints \eqref{master_c4} and $\pi_i^{\nu} \geq 0, i \in V$ for constraints \eqref{master_c5}. The reduced cost for an assignment $p$ for a (not new) daily technician $t$ is:

\begin{small}
\begin{align}
    c_p - \pi_t^{\gamma} \delta_{ot}^p - \sum_{s\in S: s\not\in S(t)} \pi_t^{\mu} \delta_{s}^p - \sum_{i \in V} \pi_i^{\nu} \delta_i^p - \Bigl( \pi_t^{\alpha} \Bigr) \label{redcost}
\end{align}
\end{small}
The reduced cost for an assignment $p$ for a new daily technician $t$ includes the dual of \eqref{master_c2} and hence becomes
\begin{small}
\begin{align}
    c_p - \pi_t^{\gamma} \delta_{ot}^p - \sum_{s\in S: s\not\in S(t)} \pi_t^{\mu} \delta_{s}^p - \sum_{i \in V} \pi_i^{\nu} \delta_i^p - \Bigl( \pi_t^{\alpha} + \pi_t^{\beta} \Bigr) \label{redcostnew}
\end{align}
\end{small}

The terms in parenthesis of the reduced costs are constants when generating new assignments. Thus the same mathematical formulation solves the assignment problem for both existing and new daily technicians. 
Introducing decision variables $u_{ot}^t, u_s^t, x_i^t$ for the parameters $\delta_{ot}^p, \delta_s^p, \delta_i^p$ 
the subproblem for a daily technician $t\in T$ with master technician $m\in \xM$ is:
\begin{small}
\begin{align}
 \min \quad & \sum_{t \in T} z_t - \pi_t^{\gamma} u_{ot} - \sum_{\begin{subarray}{c}s\in S:\\ s\not\in S(t)\end{subarray}} \pi_t^{\mu} u_s^m -  \sum_{i \in V} \pi_i^{\nu} x_i^t \rule{-5cm}{0mm}
\label{sub_obj}\\
     \text{s.t.} \quad & z_t \geq \xk \cdot c_{ij} (x_i^t + x_j^t - 1) & \forall i, j \in V \label{sub_c1} \\
      & z_t \geq \xk \cdot c_{d_ti} x_i^t & \forall i \in V \label{sub_c2}\\
    & \sum_{i \in V}f_i x_i^t \leq Q_t + w_{ot} \cdot u_{ot}^t - \sum_{i \in V} c_{d_ti}x_i^t & \label{sub_c3}\\
    & x_i^t + x_j^t \leq 1 & \forall t \in T, \forall i,j\in V: overlap(i,j) \label{sub_c3b}\\
    & x_i^t \leq  u_s^m & \forall i \in V, \forall s \in S(i): s \not\in S(t) \label{sub_c4}\\
      & x_i^t = 0 & \forall i \in V, tw(i, t) = \emptyset \label{sub_c5}\\
      & x_i^t \in \{0, 1\} & \forall i \in V \label{sub_b1}\\
      & u_{ot}^t \in \{0, 1\} & \forall a\in A, \forall t\in T\label{sub_b2}\\
      & u_s^m \in \{0, 1\} & \forall a\in A, \forall m\in \xM\label{sub_b3}
\end{align}
\end{small}
\noindent
The objective function \eqref{sub_obj} minimizes the reduced costs.
The subproblem is similar to the full investment problem \eqref{fullinv_obj} - \eqref{fullinv_b6}, but for a single daily technician, so it does not require that each task is handled, and new daily technician costs are handled implicitly. 

The subproblem is solved by a MIP solver. 
The column generation procedure is initialized by a simple heuristic, which greedily assigns tasks to each technician.
To reduce the run time of the approach, a limit is set on the number of column generation iterations. Furthermore, the number of columns generated in each iteration may be limited. Once terminated, the master problem solution may be fractional. An integer solution is solved by applying a MIP solver to the non-LP-relaxed master problem with the generated columns.

\section{Data}
\label{data}

The solution method is validated on instances from the literature \cite{mathlouthi_mip} and on an instance from the Danish Telecom infrastructure company, TDC Net.
All instances are available at Zenodo \cite{gamst2022data}.

\subsection{Data from the literature}
The benchmark instances from the literature are proposed by Mathlouthi et al. \cite{mathlouthi_mip}. They span a single day and solve a TRSP problem, where tools are taken into account. The instances thus contain information on available tools in each vehicle and on tool depots, where technicians can restock their vehicle. Also, the instances contain both travel times and travel distances. We interpret and modify the instances to match the problem considered in this paper as follows:
\begin{itemize}
    \item Tool requirements and tool depots are ignored
    \item Travel distances are ignored, since we estimate the daily travel time of a technician on the basis of the most distant tasks.
    \item Technicians work from 9 - 17
    \item Tasks have multiple time windows in the original data. We instead consider a task for each time window, so if the original data had 10 tasks with 3 time windows each, we have 30 tasks with a single time window
    \item The penalty for not serving a task is assumed to be in the same unit as the travel time
\end{itemize}
Without tool requirements and depots, groups of the Mathlouthi instances end up being identical instances. For this reason, we only consider the instances in the groups: \emph{AllSkills, nbrTech, ConfigurationDeBase, RedSkills, TpsRep10-20}. In the first group, all technicians can serve all tasks. In the remaining groups, this is not the case and the skill distribution differs between the groups. Furthermore, the number of technicians differs in \emph{nbrTech}, and task durations are shorter in \emph{TpsRep}. Finally, each group contains instances with different time window sizes (narrow and wide), and different number of tasks and technicians. An overview of the instances is provided in Table 
\ref{tab:math_num_inst_techs}. 
The total number of investigated Mathlouthi instances is 1185.

\begin{table}
\renewcommand{\tabcolsep}{0.65mm}
\scriptsize
    \centering
    \begin{tabular}{l|rrrrrrrrrrrr}\hline\hline
            &30 &45 &60 &75 &90 &105 &120 &135 &150 &210 &300 &600 \\
    \myemph{Number of inst.}& tasks& tasks& tasks& tasks& tasks& tasks& tasks& tasks& tasks& tasks& tasks& tasks\\ \hline
Narrow&50&50&50&50&50&50&50&50&90&10&50&45\\
Wide&50&50&50&50&50&50&50&50&90&10&50&40\\\hline
AllSkills&20&20&20&20&20&20&20&20&20&0&0&0\\
ConfigurationDeBase&20&20&20&20&20&20&20&20&100&20&100&85\\
nbrTech&20&20&20&20&20&20&20&20&20&0&0&0\\
RedSkills&20&20&20&20&20&20&20&20&20&0&0&0\\
TpsRep10-20&20&20&20&20&20&20&20&20&20&0&0&0\\
\hline\hline
\end{tabular}
\bigskip

\renewcommand{\tabcolsep}{1.5mm}
\scriptsize
    \centering
    \begin{tabular}{l|rrrrr}\hline\hline
    \myemph{Number of inst.}&3 tech & 4 tech & 6 tech&12 tech&24 tech\\\hline
Narrow&390&90&50&40&25\\
Wide&390&90&50&40&20\\\hline
AllSkills&180&0&0&0&0\\
ConfigurationDeBase&240&0&100&80&45\\
nbrTech&0&180&0&0&0\\
RedSkills&180&0&0&0&0\\
TpsRep10-20&180&0&0&0&0\\
\hline\hline
\end{tabular}
\caption{\label{tab:math_num_inst_techs}The number of Mathlouthi instances in the different categories.
   The first two rows report the number of instances with narrow/wide time windows. The remaining five rows report the number of instances for each instance type.}
\end{table}

The instances are useful for evaluating the performance of the ALNS and the assignment approach with and without investments. Because the instances only span a single day, the investment decisions will obviously not be realistic; however, the instances can still be used to evaluate the proposed method, especially because the instance set is large and spans different scenarios. 

\subsection{TDC Net data}
The TDC Net instance consists of 720 daily technicians, 2677 tasks, 31 different skills and a time horizon of 5 days. The instance is based on real-life data from the Danish Telecom infrastructure company, TDC Net. 
Daily technicians work on a given day from 7:30 to 15:30. They hold a number of skills; about 10\% hold less than 10 different skills, about 50 \% less than 15 skills, about 85\% less than 20 skills, and the remaining 15 \% up to 27 skills. The least frequent skill required by a task is held by 30 daily technicians, the most frequent skill by 645 daily technicians, and on average a required skill is held by 317 daily technicians.

The task time windows span from 30 minutes to all five days. About 25\% of the tasks have duration less than 30 minutes, about 33\% have duration between 30 and 60 minutes, about 30\% have duration between 60 and 120 minutes, and about 12\% have duration higher than 120 minutes. Each task requires exactly one skill. 18 of the skills are required by less than 25 tasks each, while the five most popular skills are required by 209 to 559 tasks each.

\subsection{Investment data}
\label{investment_data}
The TRSP instances from the literature and from TDC Net are extended with investment data. Recall that the objective function without investment is the sum of the total travel time and the penalty for unserved/unassigned tasks. To ensure that the objective function costs are comparable, we convert investment costs to travel time. The travel time is given in minutes. Assume that a van drives 10 km per liter of fuel and that the fuel costs \xeuro 2 per liter. Also, assume that the average driving speed is 60 km/hour. This means that \xeuro 1 equals 5 minutes of travel time.

We assume that the monthly technician salary including overhead is \xeuro 5000/month for the company. In round numbers, this gives an hourly cost of \xeuro 30, which again corresponds to 150 minutes of travel time.

In the Mathlouthi instances, the investment data is:
\begin{itemize}
    \item \myemph{Digitization}: the cost for digitizing a task is set to \xeuro 500, which corresponds to 2500 minutes of travel time. We assume that every fifth task can be digitized.
    \item \myemph{Overtime}: we assume that overtime lasts for 120 minutes. Overtime is assumed to cost 50\% extra, i.e., the total cost for the company becomes \xeuro 45 $\cdot$ 2 = \xeuro 90. This corresponds to 450 minutes of travel time.
    \item \myemph{New technician}: a technician works 8 hours per day, so the daily technician cost is \xeuro 240 which equals 1200 minutes of travel time. The new technicians consist of a duplicate of the existing technicians. E.g., if an instance consists of 2 technicians with certain skills and locations, then it is possible to hire 2 new technicians with the same skills and locations as the existing.
    \item \myemph{Skill upgrade}: assume a total of 5 workdays for training and that the training itself costs \xeuro 2000. The total cost is then \xeuro 30 $\cdot$ 8 $\cdot$ 5 + \xeuro 2000 = \xeuro 3200. We assume that the technician remains employed for another two years, giving a yearly cost of \xeuro 1600. Assuming 219 workdays per year,  the rounded daily cost is \xeuro 7. This corresponds to 35 minutes of travel time per day. 
    We generate 5 different sets of skills based on k-means clustering. 
\end{itemize}

In the TDC Net instances, the investment data is:
\begin{itemize}
    \item \myemph{Digitization}: the cost for digitizing a task is set to \xeuro 100, which corresponds to 500 minutes of travel time. We assume that tasks requiring one of skills "ADSL.I" "ADSL.S", "ISDN2.I", "ISDN2.S", "PSTN.I", "PSTN.S", "SHDSL" can be digitized, which is a total of 450 tasks.
    \item \myemph{Overtime}: the same as for the Mathlouthi instances
    \item \myemph{New technician}: the same as for the Mathlouthi instances
    \item \myemph{Skill upgrade}: the cost is the same as for the Mathlouthi instances. The daily cost is multiplied with 5 ($35\text{ minutes} \cdot 5 = 175$), because the instance spans 5 days. We generate 10 different sets of skills based on k-means clustering. 
\end{itemize}
In both sets of instances, the distance between two skills is the total number of daily technicians subtracted the number of daily technicians having both skills.

\section{Computational results}
\label{computations}

The computational evaluation is threefold:
\begin{itemize}
\item First the ALNS is benchmarked against a full MIP formulation for the TRSP on small test instances. The purpose is to ensure the quality of the ALNS.
\item First, the task assignment approach is compared to the ALNS without investment. The purpose is to certify that the assignment approach is a fair approximation of the TRSP.
\item Next, the task assignment approach with investments is evaluated. This is done by solving the corresponding ALNS, with and without investment, and comparing if the investments are actually beneficial.
\end{itemize}
In the following, we denote the assignment approach solved by column generation as the ASSM and we refer to columns in the ASSM as columns, assignments or paths.

The quality of the ALNS algorithm was assessed in the technical report \cite{gamst2022tech} by benchmarking the ALNS against a full MIP formulation for the TRSP on small test instances. The results showed that ALNS solves the TRSP with an average gap from the optimal MIP solution far below 1\%. Furthermore, the ALNS running times scale well.

The ALNS parameters were set as in \cite{gamst2022tech}: The algorithm removes between 40\% and 60\% of tasks in the destroy methods (the number of tasks to remove is chosen at random) and tries to insert all unserved task in the repair methods. The record-to-record acceptance method is run with a start threshold of 0.0015 and the decay for calculating destroy and repair weights is set to 0.99.
For the 1200 second runs, the ALNS is run with 6 retries of 200 seconds. For the 300 second runs, the ALNS is run with 3 retries of 100 seconds.

All computational evaluations are performed on a server with Intel Xeon e3-2660 CPU and 128 GB RAM. Gurobi 9.5 is used as solver. In the tables in the following sections, a ''-'' indicates that no results are available in the particular category.

All instances and full results are available at Zenodo \cite{gamst2022data}.

\subsection{Evaluation of the Adaptive Large Neighborhood Search}
\label{evaluation_alns}
The ALNS is compared to the TRSP MIP from Section~\ref{problem_definition} on the smaller Mathlouthi instances consisting of 30-75 tasks and 3-4 technicians with a time limit of 1200 seconds, see Table~\ref{tab:mip_instances}. The goal of the comparison is to quantify that the ALNS gives reasonable, albeit heuristic results and that the ALNS scales better than solving the TRSP MIP model. 
The instances are also solved by the ALNS with a time limit of 1200 seconds and of 300 seconds to test the scalability of the approach.

\begin{table}[htbp]
\scriptsize
    \centering
    \begin{tabular}{l|rrrrr}\hline\hline
        \myemph{Number of instances}&  30 tasks & 45 tasks & 60 tasks & 75 tasks & 150 tasks\\\hline
         Narrow TW&50&50&50&50&0\\
         Wide TW&50&50&50&50&5\\\hline
         AllSkills &20&20&20&20&5\\
        ConfigurationDeBase &20&20&20&20&0\\
        nbrTech &20&20&20&20&0\\
        RedSkills &20&20&20&20&0\\
        TpsRep10-20&20&20&20&20&0\\\hline\hline
\multicolumn{6}{l}{\textbf{Smaller Mathlouthi instances}}\klemme 
  \end{tabular}
    \caption{Selected instances to compare the TRSP MIP model with ALNS. The first two rows report the number of instances with narrow resp. wide time windows. The remaining five rows contain the number of instances for each instance type. Each column counts the number of instances with a given number of tasks. A total of 405 instances are selected for the comparison.}
    \label{tab:mip_instances}
\end{table}

\subsubsection{TRSP MIP results}
Results for solving the selected instances with the TRSP MIP are summarized in Table~\ref{tab:mip_results}. 
The results are reported according to the number of tasks in the instances. 
The table contains four sections. The first three sections report results according to the size of the task time windows (narrow or wide).
The last section reports results according to the instance type.
The first section holds the percentage of instances, which timed out. Then follows the average MIP gap for instances, which timed out. The last two sections contain the average time spent in seconds according to time window size resp. instance type.

\begin{figure}
    \centering
    \includegraphics[width=0.8\textwidth]{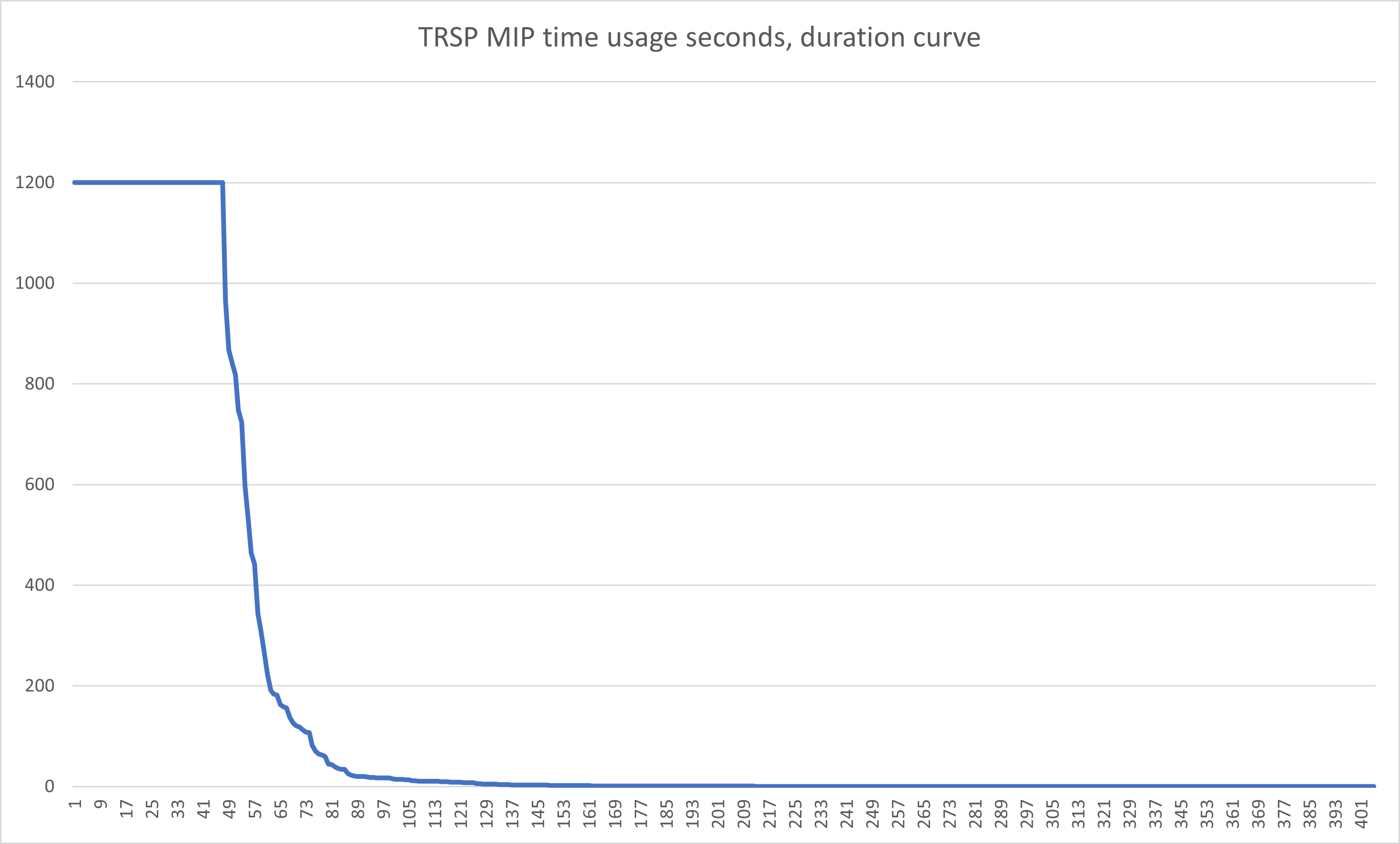}
    \caption{Duration curve for the total time spent by the TRSP MIP to solve the smaller Mathlouthi instances. Note the time limit of 1200 seconds.
             The instances (x-axis) are sorted according to nonincreasing total time.}
    \label{fig:mip_total_time_duration}
\end{figure}

The results reveal that the TRSP MIP does not scale well and that it runs out on instances containing as few as 60 tasks and on all instances with 150 tasks. The MIP gap increases dramatically from 60 tasks to 75 tasks and the time usage increases significantly when the number of tasks increase. Both indicate a significant higher complexity when the number of tasks grows larger. The results also indicate that instances with wide time windows are harder than instances with narrow windows, and that instances with many skills are also more difficult. This is most likely because the MIP solver must explore a larger solution space. 
Figure~\ref{fig:mip_total_time_duration} shows a duration curve of the total time spent. 
This again indicates that the approach scales poorly.

All in all, it is fair to conclude that the TRSP MIP approach is not feasible for larger instances. Recall that the real-life instance contains over 2500 tasks and 700 instances and that the TRSP MIP already times out on instances with 3 technicians and 150 tasks. 

\begin{table}[htbp]
\scriptsize
    \centering
    \begin{tabular}{l|rrrrr}\hline\hline
         \myemph{Pct timeouts} &  30 tasks & 45 tasks & 60 tasks & 75 tasks & 150 tasks\\\hline
         Narrow TW& 0\% & 0\% & 0\% & 6\% & -\\
         Wide TW & 0\% & 0\% & 38\% & 40\% & 100\%\\
\hline\hline
        \myemph{Avg. MIP Gap $>0$\%} &  30 tasks & 45 tasks & 60 tasks & 75 tasks & 150 tasks\\\hline
        Narrow TW & - & - & - & 8.46\% & -\\
        Wide TW &- & - & 4.86\% & 11.48\% & 28.81\%\\
\hline\hline
        \myemph{Avg. time sec.} &  30 tasks & 45 tasks & 60 tasks & 75 tasks & 150 tasks\\\hline
        Narrow TW & 0.10 & 0.66 & 4.63 & 127.68 & -\\
        Wide TW & 0.51&43.76&469.11&588.61&1200.28\\
\hline\hline
        \myemph{Avg. time sec.} &  30 tasks & 45 tasks & 60 tasks & 75 tasks & 150 tasks\\\hline
        AllSkills &0.34&48.79&564.81&605.77&1200.28 \\
        ConfigurationDeBase & 0.16&0.69&6.70&164.24& -\\
        nbrTech &  0.10&0.42&4.02&106.57& -\\
        RedSkills & 0.02&0.08&0.36&1.80 & -\\
        TpsRep10-20& 0.90&61.07&608.49&912.35 & -\\\hline\hline
\multicolumn{6}{l}{\textbf{Smaller Mathlouthi instances}}\klemme 
    \end{tabular}
    \caption{Summarized results, solving with the TRSP MIP model}
    \label{tab:mip_results}
\end{table}

\subsubsection{ALNS vs. TRSP MIP}
The proposed ALNS is compared to the TRSP MIP on the smaller instances to evaluate the performance. Both methods are run with a time limit of 1200 seconds. Summarized results can be seen in Table~\ref{tab:mip_alns_results_a} and \ref{tab:mip_alns_results_b}. 

\begin{figure}
    \centering
    \includegraphics[width=0.8\textwidth]{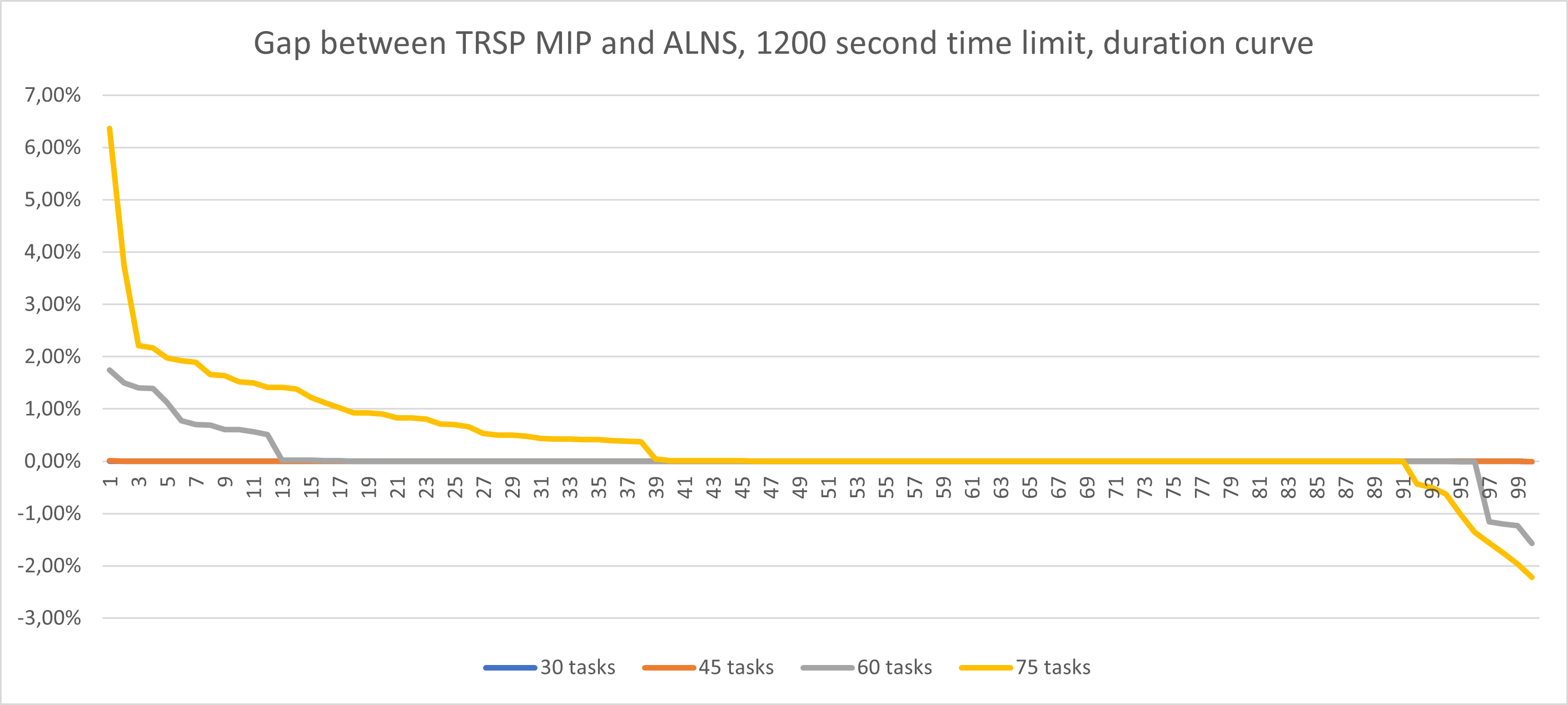}
    \caption{Duration curve of the solution value gap between the TRSP MIP and the ALNS, both with time limit of 1200 seconds, applied to the smaller Mathlouthi instances. The gap is computed as (ALNS-MIP)/MIP.}
    \label{fig:mip_alns1200_gap}
\end{figure}

The first table compares the average solution value gaps; the upper section according to the task time window sizes, and the lower section according to the instance type. A duration curve of the solution value gap for each instance is also plotted in Figure~\ref{fig:mip_alns1200_gap}. 
The results show that the ALNS finds very good solutions. The solutions are close to optimal and at times even better than the non-optimal solution of the TRSP MIP in case of timeouts. Generally, the more tasks, the larger gaps, but the gaps remain small.
The second table gives more insight into the solution differences with respect to task time window sizes. The results reveal that the average path lengths of the two approaches are very close to each other. The same goes for the average number of unserved tasks. The ALNS has a slight tendency to more unserved tasks, but not enough to cause a significant solution value gap.

\begin{table}[htbp]
\scriptsize
    \centering
    \begin{tabular}{l|rrrrr}\hline\hline
        \myemph{Avg. gap} &  30 tasks & 45 tasks & 60 tasks & 75 tasks & 150 tasks\\\hline
        Narrow TW& 0.00\%&0.00\%&0.16\%&0.51\% & -\\
        Wide TW & 0.00\%&0.00\%&-0.03\%&0.20\%&1.61\%\\
\hline\hline
        \myemph{Avg. gap} &  30 tasks & 45 tasks & 60 tasks & 75 tasks & 150 tasks\\\hline
        AllSkills & 0.00\%&0.00\%&0.11\%&0.76\%&1.61\%\\
        ConfigurationDeBase & 0.00\%&0.00\%&0.03\%&0.25\%& -\\
        nbrTech &  0.00\%&0.00\%&0.00\%&0.27\%& -\\
        RedSkills &  0.00\%&0.00\%&0.00\%&0.00\%& -\\
        TpsRep10-20&  0.00\%&0.00\%&0.19\%&0.50\%& -\\\hline\hline
\multicolumn{6}{l}{\textbf{Smaller Mathlouthi instances}. Gap $=$ (ALNS-MIP)/MIP}\klemme 
    \end{tabular}
    \caption{Gap between solution values for the TRSP MIP model and the ALNS method, both with time limit 1200 seconds}
    \label{tab:mip_alns_results_a}
\end{table}

\begin{table}[htbp]
\scriptsize
    \centering
    \begin{tabular}{l|l|rrrrr}\hline\hline
        \multicolumn{2}{l|}{\myemph{Avg. path length}}&  30 tasks & 45 tasks & 60 tasks & 75 tasks & 150 tasks\\\hline
        ALNS&Narrow TW& 5.37 & 6.55 & 7.34 & 8.06& -\\
        MIP&Narrow TW& 5.18 & 6.55 & 7.37 & 8.10& -\\
        ALNS&Wide TW & 5.12 & 6.99 & 7.60 & 8.30 & 8.73\\
        MIP&Wide TW & 4.85 & 6.67 & 7.59 & 8.26 & 8.87\\
\hline\hline
        \multicolumn{2}{l|}{\myemph{Avg. number of unserved tasks}} &  30 tasks & 45 tasks & 60 tasks & 75 tasks & 150 tasks\\\hline
        ALNS&Narrow TW&  13.60 & 24.26 & 36.74 & 49.46& -\\
        MIP&Narrow TW& 13.60 & 24.26 & 36.66 & 49.34& -\\
        ALNS&Wide TW& 14.68 & 23.88 & 35.92 & 48.70 & 123.80\\
        MIP&Wide TW & 14.68 & 23.88 & 35.92 & 48.80 & 123.40\\
\hline\hline
\multicolumn{6}{l}{\textbf{Smaller Mathlouthi instances}}
    \end{tabular}
    \caption{Comparison of path lengths and the number of unserved tasks for the TRSP MIP model and the ALNS method, both with time limit 1200 seconds.}
    \label{tab:mip_alns_results_b}
\end{table}

\begin{figure}
    \centering
    \includegraphics[width=0.8\textwidth]{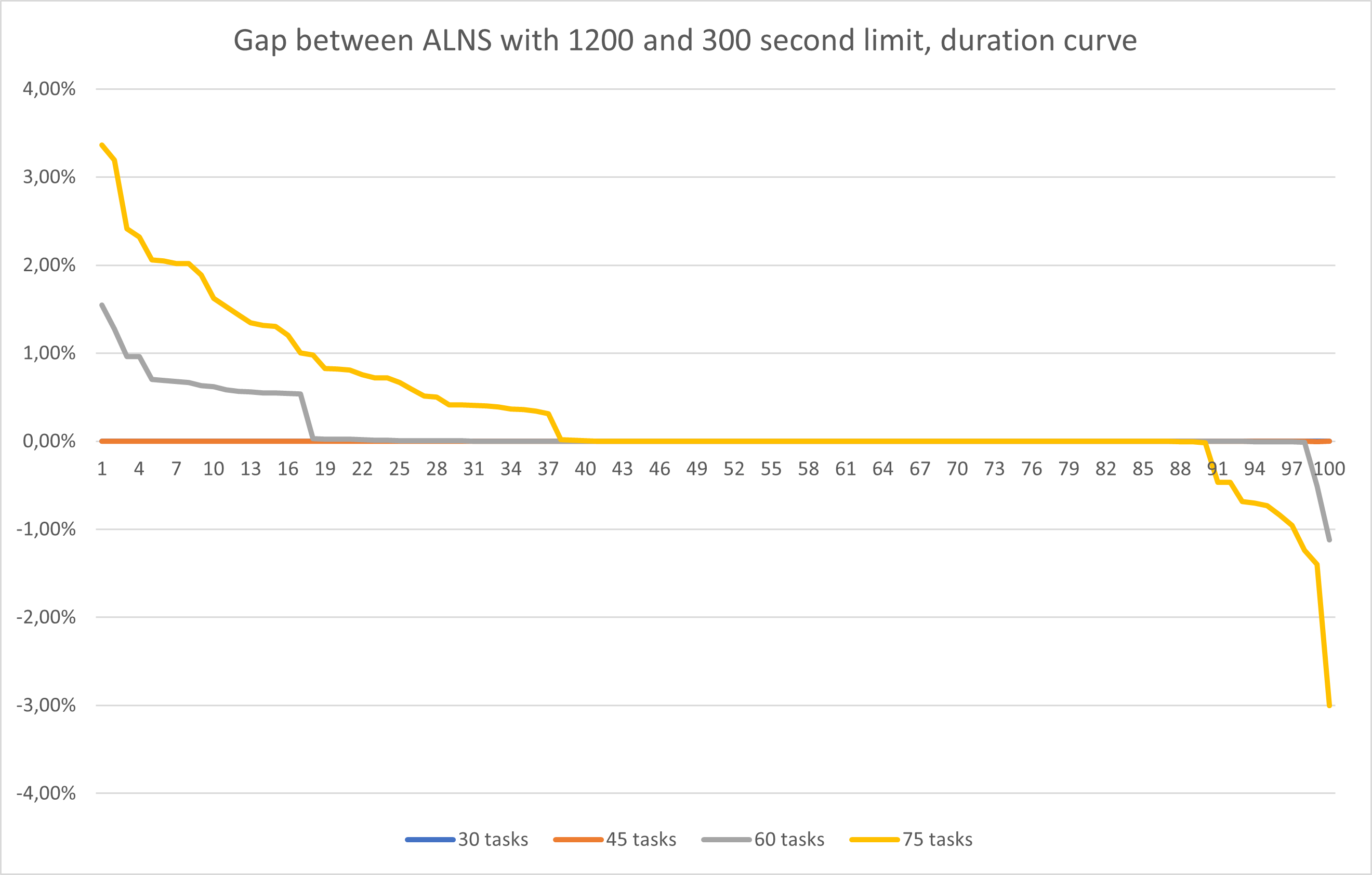}
    \caption{Duration curve of the solution value gap between the ALNS with time limit of 1200 seconds resp. 300 seconds, applied to the smaller Mathlouthi instances. The gap is computed as (ALNS300-ALNS1200)/ALNS1200.}
    \label{fig:alns1200_alns300_gap}
\end{figure}

The ALNS is also run with a time limit of 300 seconds to test its scalability. The solution value gap between the ALNS runs with 1200 seconds and 300 seconds is summarized in Table~\ref{tab:alns_alns_results}. A duration curve of the gap for each instance is also plotted in Figure~\ref{fig:alns1200_alns300_gap}. As can be seen, the average gaps are very small and the gap for each instance is only rarely more than a few percent. This indicates that the ALNS converges before 300 seconds and that the method scales well.

All in all, the gaps reveal that the ALNS solves the TRSP very well and scales well.
 We thus conclude that comparing the ASSM with and without investments to results from the ALNS with a 300 second time limit, is a fair approach to evaluate performance.

\begin{table}[htbp]
\scriptsize
    \centering
    \begin{tabular}{l|rrrrr}\hline\hline
        \myemph{Avg. gap} (in \%) &  30 tasks & 45 tasks & 60 tasks & 75 tasks & 150 tasks\\\hline
        Narrow TW& 0.00 & 0.00 & 0.15 & 0.33& -\\
        Wide TW & 0.00 & 0.00 & 0.08 & 0.33 & 0.35\\
\hline\hline
        \myemph{Avg. gap} (in \%) &  30 tasks & 45 tasks & 60 tasks & 75 tasks & 150 tasks\\\hline
        AllSkills & 0.00 & 0.00 & 0.21 & 0.30 & 0.35\\
        ConfigurationDeBase & 0.00 & 0.00 & 0.00 & 0.33& -\\
        nbrTech & 0.00 & 0.00 & 0.08 & 0.10& -\\
        RedSkills & 0.00 & 0.00 & 0.00 & 0.02& -\\
        TpsRep10-20& 0.00 & 0.00 & 0.26 & 0.90& -\\\hline\hline
\multicolumn{6}{l}{\textbf{Smaller Mathlouthi instances}. Gap $=$ (ALNS300-ALNS1200)/ALNS1200}\klemme 
    \end{tabular}
    \caption{Gap between solution values for the ALNS method with time limit 1200 seconds resp. 300 seconds.}
    \label{tab:alns_alns_results}
\end{table}

\subsection{Evaluation of the ASSM}

The results from the ASSM are compared to that of the ALNS. The comparison is performed on all instances. We compare the objective function values, the number of unserved tasks, path length and number of used technicians. The costs in the objective functions are expected to differ; the ALNS optimizes the actual travel time while the ASSM optimizes an estimated 
travel time. Still, it is relevant to see how close the two values are, because investments depend on the estimated travel time in the ASSM approach.

Both the ASSM and the ALNS are given a time limit of 300 seconds.
Also, we reduce the set of tasks, which can be assigned to a technician, by setting an upper bound of 100 minutes in travel time from the technician depot to the task location. This speeds up both approaches, as fewer tasks need to be considered for each technician. 

In the following, we first analyze the results on the Mathlouthi instances, then on the TDC Net instance.

\subsubsection{Results for the Mathlouthi instances}
\label{mathlouthi_bc_results}
The ASSM parameters are set as follows for the Mathlouthi instances: The scaling factor $\xk$ is set to 5 as this is a good estimate of the additional driving between the outer-most tasks. The number of column generation iterations is limited to 75 and in each iteration at most one path is generated per technician. The subproblem is given a 3 second time limit and the time limit for solving the master problem in each iteration is 30 seconds. At most 60 seconds can be spent on solving the integer master problem at the very end of the procedure.

The ALNS is tuned as explained in the beginning of 
Section~\ref{computations}
with 3 retries of 100 seconds. 

\begin{figure}
\centering
    \includegraphics[width=0.8\textwidth]{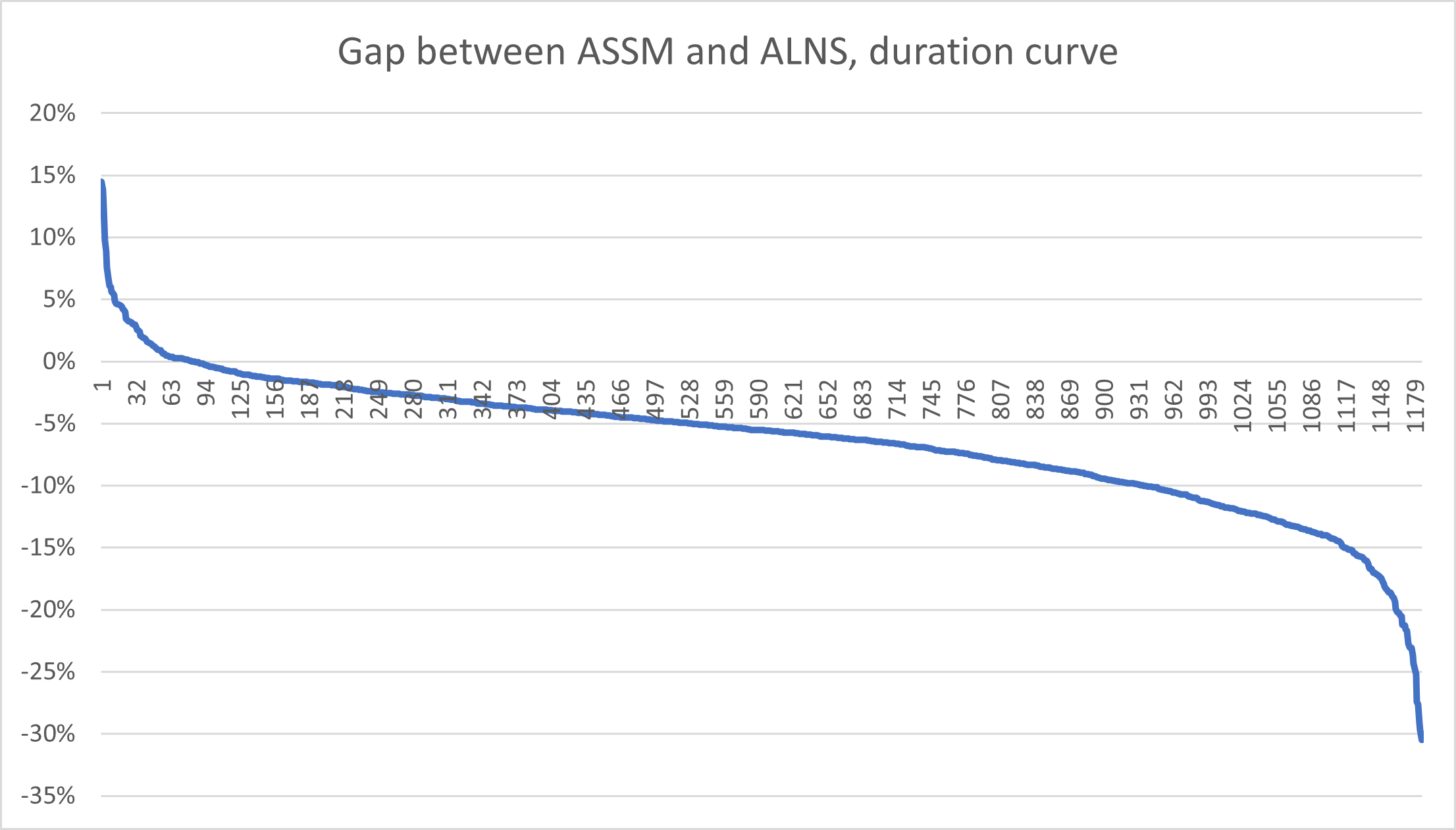}
    \caption{Duration curve of the solution value gap between the ASSM and the ALNS, applied to the Mathlouthi instances. The gap is computed as (ASSM-ALNS)/ALNS.
             The instances (x-axis) are sorted according to nonincreasing gaps.}
    \label{fig:assm_alns_gap}
 \end{figure}

Tables \ref{tab:alns_assm_results_tasks} and \ref{tab:alns_assm_results_techs} show the average gap between the ASSM and the ALNS solution values, according to the number of tasks resp. number of technicians. The gaps indicate that the ASSM finds smaller solution values than the ALNS.  The gaps do not seem to increase with the size of the instances, but they are generally larger on the instances with narrow time windows. This pattern was also seen when comparing the solution values of ALNS and TRSP MIP: Narrow time windows restrict the solution space, which benefits the MIP solver and challenge ALNS. Thus the ASSM on average finds much lower solution values for these instances than the ALNS.
Still, the duration curve of solution value gaps in Figure~\ref{fig:assm_alns_gap} shows that the gap of single instances are within reasonable bounds of between appr. -30\% and 15\%. 

Tables \ref{tab:alns_assm_unserved_tasks} and \ref{tab:alns_assm_unserved_techs} report the relative difference between number of unserved tasks, according to the number of tasks resp. number of technicians. The negative values mean that the ALNS on average has more unserved tasks than the ASSM. This is probably the main reason for the difference in solution values. There is no clear pattern between problem instance size and the difference.

Tables \ref{tab:alns_assm_path_tasks} and \ref{tab:alns_assm_path_techs} contain the relative difference between the path lengths, according to the number of tasks resp. number of technicians. The ASSM generally finds longer paths than the ALNS, probably because it generates less accurate schedules and thus may produce infeasible routes. This explains why it is capable of serving more tasks and has better solution values than the ALNS.

Finally, Tables \ref{tab:alns_assm_tech_tasks} and \ref{tab:alns_assm_tech_techs} contain the relative difference between the number of used daily technicians, according to the number of tasks resp. number of technicians. The two methods agree on the number of used daily technicians to a large extend. 

All in all, the results show that the ASSM finds solutions, which compare reasonably well to the ALNS. 
The mathematical formulation for the subproblem constrains the length of the paths by estimating the time spent on travelling (constraints \eqref{sub_c3}) and that two overlapping tasks cannot be serviced by the same technician (constraints \eqref{sub_c4}). Preliminary tests showed no major difference in the results if the estimated travel time was increased significantly. Instead the results point towards missing information about three or more tasks overlapping. For example, let tasks $a$, $b$ and $c$ all have time windows $[8:12]$ on the same day and have duration $f_a = f_b = f_c = 2$. Constraints \eqref{sub_c4} only check for pair of tasks and the mathematical formulation for the subproblem will thus not detect that these three tasks cannot be serviced by the same technician. Thus the ASSM ends up with assignments that are infeasible in the ALNS and the ASSM results are too optimistic compared to the ALNS.
The ASSM subproblem formulation could be extended to consider groups of three or more tasks, either by enumerating conflicting tasks up front or through a cutting plane procedure. Both approaches would, however, increase the time spent on solving the subproblem.


\begin{table}
\renewcommand{\tabcolsep}{0.4mm}
\scriptsize
    \centering
    \begin{tabular}{l|rrrrrrrrrrrr}\hline\hline
            &30 &45 &60 &75 &90 &105 &120 &135 &150 &210 &300 &600 \\
\myemph{Avg obj gap} (in \%)& tasks& tasks& tasks& tasks& tasks& tasks& tasks& tasks& tasks& tasks& tasks& tasks\\\hline
Narrow&-8.09&-12.26&-10.40&-8.82&-8.13&-9.64&-7.49&-6.96&-7.05&-2.10&-8.78&-6.67\\
Wide&-3.14&-4.00&-3.53&-3.51&-3.18&-4.28&-4.70&-5.09&-4.37&-2.01&-6.40&-4.58\\
\hline
AllSkills&-4.93&-12.09&-12.78&-8.53&-7.54&-6.67&-4.84&-4.13&-4.13&-&-&-\\
ConfigurationDeBase&-8.01&-8.92&-6.24&-6.48&-5.11&-5.30&-4.40&-4.32&-5.21&-2.06&-7.59&-5.69\\
nbrTech&-7.20&-10.43&-9.61&-7.55&-6.82&-8.51&-7.31&-6.08&-5.14&-&-&-\\
RedSkills&-5.55&-8.65&-6.39&-4.72&-3.40&-4.20&-2.98&-3.16&-3.14&-&-&-\\
TpsRep10-20&-2.39&-0.56&0.19&-3.56&-5.40&-10.13&-10.95&-12.43&-12.95&-&-&-\\
\hline\hline
\multicolumn{10}{l}{\textbf{Mathlouthi instances}. Gap $=$ (ASSM-ALNS)/ALNS}\klemme 
    \end{tabular}
    \caption{\label{tab:alns_assm_results_tasks}
    Gap between solution values for the ALNS method and the ASSM, grouped by the number of tasks}
\end{table}

\begin{table}
\renewcommand{\tabcolsep}{1.5mm}
\scriptsize
    \centering
    \begin{tabular}{l|rrrrr}\hline\hline
\myemph{Avg obj gap} (in \%)&3 tech&4 tech&6 tech&12 tech&24 tech\\\hline
Narrow&-7.92&-10.38&-8.01&-10.00&-6.48\\
Wide&-3.66&-4.88&-5.19&-7.25&-3.61\\\hline
AllSkills&-7.29&-&-&-&-\\
ConfigurationDeBase&-4.98&-&-6.60&-8.63&-5.20\\
nbrTech&-&-7.63&-&-&-\\
RedSkills&-4.69&-&-&-&-\\
TpsRep10-20&-6.46&-&-&-&-\\
\hline\hline
\multicolumn{6}{l}{\textbf{Mathlouthi instances}. Gap $=$ (ASSM-ALNS)/ALNS}\klemme 
    \end{tabular}
    \caption{\label{tab:alns_assm_results_techs}Gap between solution values for the ALNS method and the ASSM, grouped by the number of technicians. }
\end{table}

\begin{table}
\renewcommand{\tabcolsep}{0.4mm}
\scriptsize
    \centering
    \begin{tabular}{l|rrrrrrrrrrrr}\hline\hline
\myemph{Avg. number un-} &30 &45 &60 &75 &90 &105 &120 &135 &150 &210 &300 &600 \\
\myemph{served tasks} (in \%)& tasks & tasks& tasks& tasks& tasks& tasks& tasks& tasks& tasks& tasks& tasks& tasks\\\hline
Narrow&-10.80&-13.81&-9.25&-9.13&-10.21&-11.65&-9.96&-9.36&-7.85&-2.43&-5.48&-6.50\\
Wide&-4.17&-5.81&-4.30&-5.76&-4.41&-6.56&-6.89&-7.60&-5.77&-2.21&-3.80&-3.55\\\hline
AllSkills&-10.31&-19.09&-16.13&-12.55&-11.79&-8.23&-7.87&-6.13&-5.36&-&-&-\\
ConfigurationDeBase&-10.88&-12.21&-6.32&-7.93&-5.60&-6.71&-6.06&-6.37&-5.60&-2.32&-4.64&-5.11\\
nbrTech&-10.95&-12.56&-7.34&-4.22&-3.50&-6.69&-5.65&-5.27&-3.97&-&-&-\\
RedSkills&-5.33&-7.16&-3.72&-1.16&-1.99&-4.44&-4.40&-4.28&-5.06&-&-&-\\
TpsRep10-20&0.03&1.96&-0.37&-11.37&-13.65&-19.44&-18.12&-20.35&-18.91&-&-&-\\
\hline\hline
\multicolumn{10}{l}{\textbf{Mathlouthi instances}. Relative difference $=$ (ASSM-ALNS)/ALNS}\klemme 
    \end{tabular}
    \caption{\label{tab:alns_assm_unserved_tasks}
    Relative difference in unserved tasks for the ALNS method and the ASSM, grouped by the number of tasks.}
\end{table}

\begin{table}
\renewcommand{\tabcolsep}{1.5mm}
\scriptsize
    \centering
    \begin{tabular}{l|rrrrr}\hline\hline
    \myemph{Avg. number userved tasks} (in \%)&3 tech&4 tech&6 tech&12 tech&24 tech\\ \hline
Narrow&-10.09&-9.23&-5.05&-6.13&-9.70\\
Wide&-5.96&-4.14&-3.42&-4.18&-5.72\\\hline
AllSkills&-10.83&-&-&-&-\\
ConfigurationDeBase&-6.47&-&-4.23&-5.15&-7.93\\
nbrTech&-&-6.69&-&-&-\\
RedSkills&-4.17&-&-&-&-\\
TpsRep10-20&-11.14&-&-&-&-\\
\hline\hline
\multicolumn{6}{l}{\textbf{Mathlouthi instances}. Relative difference $=$ (ASSM-ALNS)/ALNS}\klemme 
    \end{tabular}
    \caption{\label{tab:alns_assm_unserved_techs}Relative difference in unserved tasks for the ALNS method and the ASSM, grouped by the number of technicians. }
\end{table}

\begin{table}
\renewcommand{\tabcolsep}{0.65mm}
\scriptsize
    \centering
    \begin{tabular}{l|rrrrrrrrrrrr}\hline\hline
\myemph{Avg. path length}&30 &45 &60 &75 &90 &105 &120 &135 &150 &210 &300 &600 \\
(in \%)         & tasks& tasks& tasks& tasks& tasks& tasks& tasks& tasks& tasks& tasks& tasks& tasks\\\hline
Narrow&22.11&21.43&17.51&19.45&24.63&32.10&31.32&31.66&28.58&18.40&16.18&18.48 \\
Wide&21.62&8.38&8.71&12.27&10.80&16.54&20.40&24.33&19.65&16.88&10.73&11.89\\\hline
AllSkills&35.15&24.47&26.88&26.86&31.33&26.89&29.66&26.10&25.46&-&-&-\\
ConfigurationDeBase&30.74&22.99&14.82&21.09&18.58&25.60&26.03&30.49&20.72&17.64&13.45&15.38\\
nbrTech&18.94&15.15&11.58&8.12&8.45&17.78&17.34&18.02&14.90&-&-&-\\
RedSkills&9.87&12.89&8.77&3.64&7.62&17.19&19.45&21.14&27.34&-&-&-\\
TpsRep10-20&14.62&-0.97&3.52&19.58&22.60&34.15&36.82&44.24&45.74&-&-&-\\
\hline\hline
\multicolumn{10}{l}{\textbf{Mathlouthi instances}. Relative difference $=$ (ASSM-ALNS)/ALNS}\klemme 
    \end{tabular}
    \caption{\label{tab:alns_assm_path_tasks}
    Relative difference in path lengths for the ALNS method and the ASSM, grouped by the number of tasks.}
\end{table}

\begin{table}
\renewcommand{\tabcolsep}{1.5mm}
\scriptsize
    \centering
    \begin{tabular}{l|rrrrr}\hline\hline
\myemph{Avg. path length} (in \%)          &3 tech&4 tech&6 tech&12 tech&24 tech\\ \hline
Narrow&27.41&18.72&15.60&17.36&22.34\\
Wide&17.73&10.23&9.90&11.57&15.36\\\hline
AllSkills&28.09&-&-&-&-\\
ConfigurationDeBase&23.28&-&12.75&14.47&19.24\\
nbrTech&-&14.47&-&-&-\\
RedSkills&14.21&-&-&-&-\\
TpsRep10-20&24.48&-&-&-&-\\\hline\hline
\multicolumn{6}{l}{\textbf{Mathlouthi instances}. Relative difference $=$ (ASSM-ALNS)/ALNS}\klemme 
    \end{tabular}
    \caption{\label{tab:alns_assm_path_techs}Relative difference in path lengths for the ALNS method and the ASSM, grouped by the number of technicians.}
\end{table}

\begin{table}
\renewcommand{\tabcolsep}{0.7mm}
\scriptsize
    \centering
    \begin{tabular}{l|rrrrrrrrrrrr}\hline\hline
\myemph{Avg. number used} &30 &45 &60 &75 &90 &105 &120 &135 &150 &210 &300 &600 \\
\myemph{technicians} (in \%) & tasks& tasks& tasks& tasks& tasks& tasks& tasks& tasks& tasks& tasks& tasks& tasks\\\hline
Narrow&-2,33&-0,67&0,00&0,00&0,00&0,00&0,00&0,00&0,00&0,00&0,00&0,00\\
Wide&-6,83&4,00&0,00&0,00&0,00&0,00&0,00&0,00&0,00&0,00&0,00&0,00\\\hline
AllSkills&-10,00&0,00&0,00&0,00&0,00&0,00&0,00&0,00&0,00&-&-&-\\
ConfigurationDeBase&-6,67&0,00&0,00&0,00&0,00&0,00&0,00&0,00&0,00&0,00&0,00&0,00\\
nbrTech&-3,75&0,00&0,00&0,00&0,00&0,00&0,00&0,00&0,00&-&-&-\\
RedSkills&0,00&0,00&0,00&0,00&0,00&0,00&0,00&0,00&0,00&-&-&-\\
TpsRep10-20&-2,50&8,33&0,00&0,00&0,00&0,00&0,00&0,00&0,00&-&-&-\\
\hline\hline
\multicolumn{10}{l}{\textbf{Mathlouthi instances}. Relative difference $=$ (ASSM-ALNS)/ALNS}\klemme 
    \end{tabular}
    \caption{\label{tab:alns_assm_tech_tasks}Relative difference in the number of used daily technicians for the ALNS method and the ASSM, grouped by the number of tasks.}
\end{table}

\begin{table}
\renewcommand{\tabcolsep}{1.5mm}
\scriptsize
    \centering
    \begin{tabular}{l|rrrrr}\hline\hline
\myemph{Avg. number used technicians} (in \%)&3 tech&4 tech&6 tech&12 tech&24 tech\\ \hline
Narrow&-0.38&0.00&0.00&0.00&0.00\\
Wide&-0.17&-0.83&0.00&0.00&0.00\\\hline
AllSkills&-1.11&-&-&-&-\\
ConfigurationDeBase&-1.11&-&0.00&0.00&0.00\\
nbrTech&-&-0.42&-&-&-\\
RedSkills&0.00&-&-&-&-\\
TpsRep10-20&0.65&-&-&-&-\\
\hline\hline
\multicolumn{6}{l}{\textbf{Mathlouthi instances}. Relative difference $=$ (ASSM-ALNS)/ALNS}\klemme 
    \end{tabular}
    \caption{\label{tab:alns_assm_tech_techs}Relative difference in the number of used daily technicians for the ALNS method and the ASSM, grouped by the number of technicians.}
\end{table}


\subsubsection{Results for the TDC Net instance}\label{results_tdc_bc}

The TDC Net instance is significantly larger than the Mathlouthi instances, so the time limit is increased to three hours. For both the ASSM and the ALNS, we keep a maximum travel time of 100 minutes from the depot to any task. In the ASSM, we also keep $k=5$. In each column generation iteration, at most 500 columns can be generated. We set this limit to promote more column generation iterations before timing out, in hopes of achieving better dual variable values and thus better columns in the latter column generation iterations.  As the number of technicians exceeds the number of generated paths per iteration, we sort the technicians according to the decreased order of the term in parenthesis of the reduced cost in constraint \eqref{redcost} (or constraint \eqref{redcostnew} for new technicians). 
Preliminary tests revealed that if we try to generate a column for each technician in each column generation iteration, the instance times out early and the number of generated columns become untractable.

We increase the time limit for solving the relaxed master problem in each column generation iteration to 2 minutes. We preserve a 3 second time limit for solving each subproblem, because the travel limit of 100 minutes from the depot limits the number of tasks to consider for each technician and thus also the subproblem size. Once the column generation procedure finishes, the integer problem is given a time limit of 20 minutes. 

The ALNS is also given more time; a total of 2 restarts each with a time limit of 5400 seconds. Preliminary testing suggests that the ALNS needs more time to search the neighborhood, hence we have reduced the number of restarts to provide more time in each iteration. Furthermore, recall that the ALNS finishes with solving a set cover like model with up to 25000 of the most promising generated path. This will further increase the total run time.

All in all, this makes both the ASSM and the ALNS rather time consuming. This should, however, not be a problem in practice, as they are used for making strategic decisions and not operational plans. 

Results for both approaches are summarized in Table~\ref{tab:tdc_no_investment_results}. They show that the ASSM times out both in the column generation and when solving the final integer master problem. Still, the solution value of the ASSM is better than that of the ALNS, because it has fewer unserved tasks. This corresponds to what we saw for the Mathlouthi instances and is probably caused by the ASSM approach being more optimistic and not realizing that groups of three or more tasks overlap. Also, the last ALNS improvements happen rather late which suggests that the ALNS also times out before convergence.

If we subtract the penalty of unserved tasks from the objective, then the travel time is significantly smaller for the ALNS. This indicates that the travel time estimate in the ASSM is overestimated. However, looking at the total objective the two solution approaches only differ by 13\%. Since both approaches time out, it is difficult to assess whether the difference stems from the different solution approaches and objectives or from lack of time.

It should be noted that in real-life, TDC Net does not leave tasks unserved. Instead, overtime would be appointed to technicians, or technicians from a neighboring area would help with the final tasks. As a final resort, a task could be pushed to the next day and given very high priority.

\begin{table}[p]
\scriptsize
    \centering
    \begin{tabular}{l|rr}\hline\hline
\myemph{Comparison}
        & ASSM & ALNS\\\hline
        Objective & 383756&441858\\
        Total time sec & 11201.78&13478.58\\
        \# unserved tasks & 16& 42\\
        Total penalty & 142836 &406357\\
        Avg. path length & 5.86 &4.70\\
        \# paths in solution & 545 &561 \\
        Total number generated paths & 19484 & 292026\\\hline
        Number cg. iters. & 38 & - \\
        Optimal cg. & False & - \\
        Optimal int. master problem & False & - \\ \hline
        Avg. number improvements & - & 4.5 \\
        Avg. number iters. & - & 1500  \\
        Avg. last impr. iter. & - & 1336 \\ \hline \hline
\multicolumn{3}{l}{\textbf{TDC Net instances}}\klemme 
  \end{tabular}
    \caption{Results for solving the
instances with the ASSM and ALNS without investments. The last two sections concern results for the ASSM resp. the ALNS only.}
    \label{tab:tdc_no_investment_results}
\end{table}

\subsection{Evaluation of investments}
The Mathlouthi and TDC Net instances are solved by the ASSM with investment options. The resulting investments are included in the instances and they are then resolved by the ALNS. This can then be compared to the ALNS solutions without investments to analyze if the approach actually produces positive business cases and sane investment suggestions. In the following, we first analyze results for the Mathlouthi instances, then for the TDC Net instance.

\subsubsection{Mathlouthi instances}
The algorithms in this section are run with the same parameter settings as explained in Section~\ref{mathlouthi_bc_results}.

The Mathlouthi instances are first solved with the ASSM including investments. Tables \ref{tab:math_inv_opt} and \ref{tab:math_inv_opt_techs} show that many of the smaller instances and few of the larger instances are solved to optimality by the column generation approach. The number of column generation iterations in tables \ref{tab:math_inv_cg} and \ref{tab:math_inv_cg_techs} decrease as the instances grow larger, due to timing out. Finally, we consider the number of generated paths in tables \ref{tab:math_inv_gen_path} and \ref{tab:math_inv_gen_path_techs}. More paths are generated for the larger instances, but not significantly more because the number of column generation iterations also decreases. The non-optimal solutions could mean that the investments decisions could be improved further.

\begin{table}
\renewcommand{\tabcolsep}{0.5mm}
\scriptsize
    \centering
    \begin{tabular}{l|rrrrrrrrrrrr}\hline\hline
\myemph{Optimal solution}&30 &45 &60 &75 &90  &105 &120 &135 &150 &210 &300 &600 \\
(in \%)         & tasks& tasks& tasks& tasks& tasks & tasks& tasks& tasks& tasks& tasks& tasks& tasks\\\hline
Narrow&96.00&82.00&76.00&70.00&70.00&54.00&34.00&34.00&18.89&0.00&0.00&0.00\\
Wide&96.00&60.00&32.00&22.00&4.00&4.00&4.00&0.00&1.11&0.00&0.00&0.00\\\hline
AllSkills&95.00&70.00&50.00&55.00&35.00&30.00&25.00&30.00&20.00&-&-&-\\
ConfigurationDeBase&100.00&70.00&45.00&30.00&35.00&35.00&10.00&5.00&5.00&0.00&0.00&0.00\\
nbrTech&90.00&55.00&25.00&30.00&15.00&15.00&5.00&0.00&5.00&-&-&-\\
RedSkills&100.00&60.00&50.00&30.00&45.00&35.00&10.00&5.00&5.00&-&-&-\\
TpsRep10-20&95.00&100.00&100.00&85.00&55.00&30.00&45.00&45.00&35.00&-&-&-\\
\hline\hline
\multicolumn{6}{l}{\textbf{Mathlouthi instances}}\klemme 
    \end{tabular}
    \caption{\label{tab:math_inv_opt}
    Percentage of instances solved to optimality by the ASSM with investments.}
\end{table}

\begin{table}
\renewcommand{\tabcolsep}{1.5mm}
\scriptsize
    \centering
    \begin{tabular}{l|rrrrr}\hline\hline
    \myemph{Optimal solution} (in \%)&3 tech & 4 tech & 6 tech&12 tech&24 tech\\\hline
Narrow &61.28&37.78&4.00&0.00&0.00 \\
Wide&24.87&15.56&2.00&0.00&0.00\\\hline
AllSkills&45.56&-&-&-&-\\
ConfigurationDeBase&28.33&-&3.00&0.00&0.00\\
nbrTech&-&26.67&-&-&-\\
RedSkills&37.78&-&-&-&-\\
TpsRep10-20&65.56&-&-&-&-\\
\hline\hline
\multicolumn{6}{l}{\textbf{Mathlouthi instances}}\klemme 
    \end{tabular}
    \caption{\label{tab:math_inv_opt_techs}Percentage of instances solved to optimality by the ASSM with investments.}
\end{table}

\begin{table}
\renewcommand{\tabcolsep}{0.6mm}
\scriptsize
    \centering
    \begin{tabular}{l|rrrrrrrrrrrr}\hline\hline
                     &30 &45 &60 &75 &90  &105 &120 &135 &150 &210 &300 &600 \\
\myemph{Avg num cg iterations}& tasks& tasks& tasks& tasks& tasks & tasks& tasks& tasks& tasks& tasks& tasks& tasks\\\hline
Narrow&32.00&38.44&42.08&48.04&54.68&62.24&66.30&69.38&72.49&75.00&24.04&2.87\\
Wide&44.24&58.12&67.34&71.40&74.76&74.80&74.48&75.00&74.87&75.00&31.80&3.95\\\hline
AllSkills&37.05&50.20&51.20&58.50&64.60&66.50&67.60&70.25&72.25&-&-&-\\
ConfigurationDeBase&37.80&53.95&62.20&63.30&65.10&68.00&73.15&74.80&74.46&75.00&27.92&3.38\\
nbrTech&48.45&55.00&64.20&66.95&72.10&72.60&74.75&75.00&74.65&-&-&-\\
RedSkills&40.45&54.45&61.10&63.60&65.70&69.05&73.90&75.00&75.00&-&-&-\\
TpsRep10-20&26.85&27.80&34.85&46.25&56.10&66.45&62.55&65.90&68.90&-&-&-\\
\hline\hline
\multicolumn{6}{l}{\textbf{Mathlouthi instances}}\klemme 
    \end{tabular}
    \caption{\label{tab:math_inv_cg}
    \mbox{Avg. number of col. generation iterations before optimality or timeout}}
\end{table}

\begin{table}
\renewcommand{\tabcolsep}{1.5mm}
\scriptsize
    \centering
    \begin{tabular}{l|rrrrr}\hline\hline
    \myemph{Avg num cg iterations}&3 tech & 4 tech & 6 tech&12 tech&24 tech\\\hline
Narrow&53.29&62.87&47.72&9.63&2.00\\
Wide&68.18&71.29&53.36&12.73&2.95\\ \hline
AllSkills&59.79&-&-&-&-\\
ConfigurationDeBase&66.39&-&50.54&11.18&2.42\\
nbrTech&-&67.08&-&-&-\\
RedSkills&64.25&-&-&-&-\\
TpsRep10-20&50.63&-&-&-&-\\
\hline\hline
\multicolumn{6}{l}{\textbf{Mathlouthi instances}}\klemme 
    \end{tabular}
    \caption{\label{tab:math_inv_cg_techs}\mbox{Avg. number of col. generation iterations before optimality or timeout}}
\end{table}

\begin{table}
\renewcommand{\tabcolsep}{0.4mm}
\scriptsize
    \centering
    \begin{tabular}{l|rrrrrrrrrrrr}\hline\hline
\myemph{Avg num}        &30 &45 &60 &75 &90  &105 &120 &135 &150 &210 &300 &600 \\
\myemph{generated paths}& tasks& tasks& tasks& tasks& tasks & tasks& tasks& tasks& tasks& tasks& tasks& tasks\\\hline
Narrow&125.36&188.06&201.44&218.72&235.98&284.50&328.70&345.66&419.99&355.20&366.68&132.07\\
Wide&121.88&230.82&321.08&363.80&410.30&402.72&395.94&405.18&480.10&400.10&483.98&155.73\\\hline
AllSkills&120.40&220.90&247.25&282.50&305.75&323.05&334.15&342.95&353.05&-&-&-\\
ConfigurationDeBase&129.30&234.45&299.55&311.15&313.85&337.50&373.20&379.25&510.75&377.65&425.33&143.20\\
nbrTech&151&256.95&331.15&374.10&429.75&448.45&474.80&490.40&492.30&-&-&-\\
RedSkills&136.55&237.80&300.30&317.25&315.90&337.45&373.95&381.35&374.15&-&-&-\\
TpsRep10-20&80.85&97.10&128.05&171.30&250.45&271.60&255.50&283.15&277.15&-&-&-\\
\hline\hline
\multicolumn{6}{l}{\textbf{Mathlouthi instances}}\klemme 
    \end{tabular}
    \caption{\label{tab:math_inv_gen_path}
    Avg. number of generated paths in col. generation}
\end{table}

\begin{table}
\renewcommand{\tabcolsep}{1.5mm}
\scriptsize
    \centering
    \begin{tabular}{l|rrrrr}\hline\hline
    \myemph{Avg num generated paths}&3 tech & 4 tech & 6 tech&12 tech&24 tech\\\hline
Narrow&240.46&333.46&490.74&253.85&142.68\\
Wide&321.90&432.97&587.46&323.95&169.30\\ \hline
AllSkills&281.11&-&-&-&-\\
ConfigurationDeBase&320.51&-&539.10&288.90&154.51\\
nbrTech&-&383.21&-&-&-\\
RedSkills308.30&-&-&-&-\\
TpsRep10-20&201.68&-&-&-&-\\
\hline\hline
\multicolumn{6}{l}{\textbf{Mathlouthi instances}}\klemme 
    \end{tabular}
    \caption{\label{tab:math_inv_gen_path_techs}Avg. number of generated paths in col. generation}
\end{table}

\bigskip

Next, we consider the actual investments in the Mathlouthi instances which are summarized in tables \ref{tab:math_inv_skill}-\ref{tab:math_inv_over_techs}. Generally, the number of investments increases when the number of tasks and technicians increases. Skill investments are generally attractive in all instances, except the \emph{AllSkills} instances where each technician already possess all skills. The highest number of skill investments are observed in the instance set \emph{RedSkills}, where the skills per technician are quite small.

The investments clearly show that many technicians are assigned overtime. Still, investments are made in new technicians, which indicate that the workload is high relative to the number of technicians. This is also supported by the number of tasks being digitized.

\begin{figure}
\centering
    \includegraphics[width=0.8\textwidth]{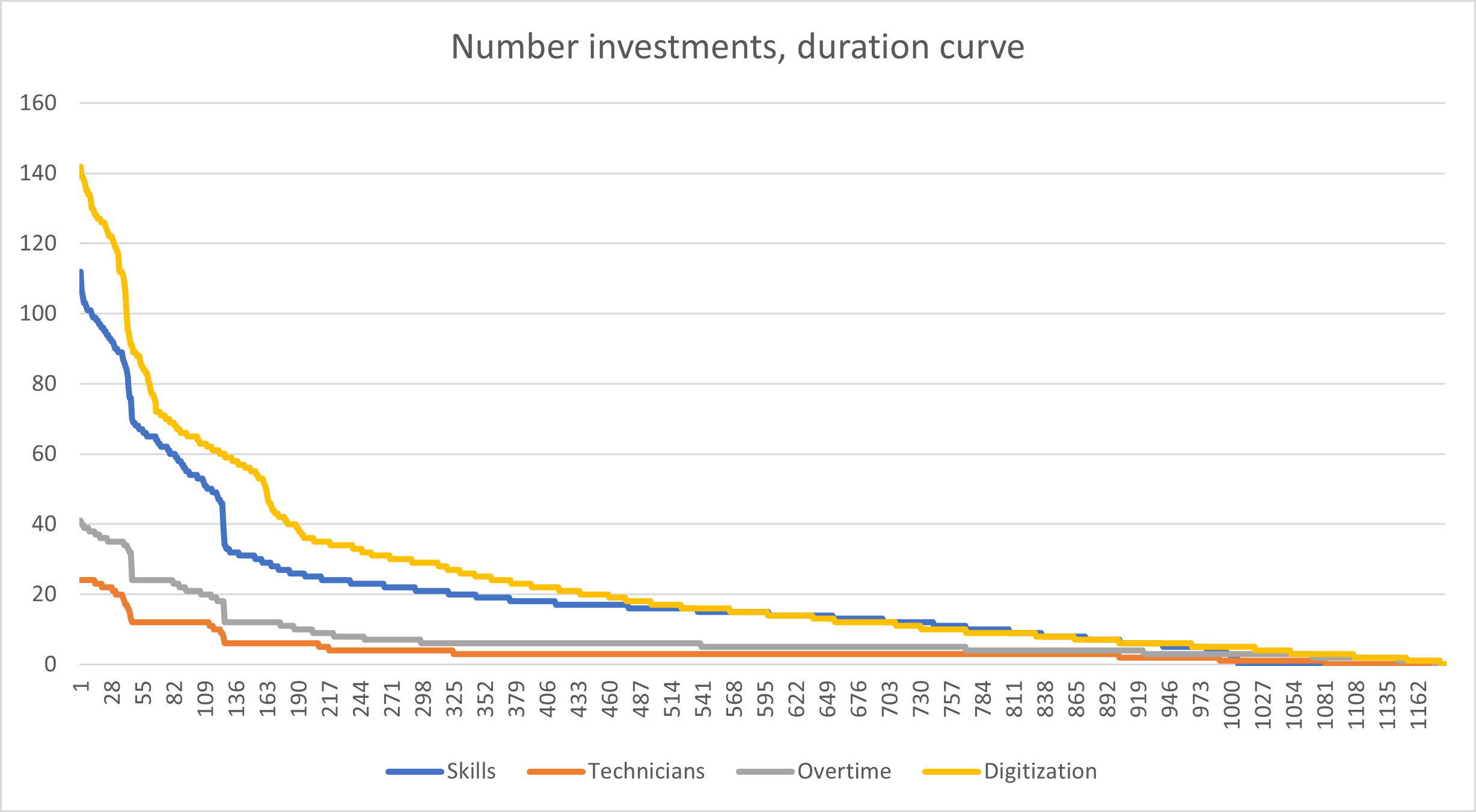}
    \caption{Duration curve of the number of investments.}
    \label{fig:num_inv}
 \end{figure}

\begin{figure}
\centering
    \includegraphics[width=0.8\textwidth]{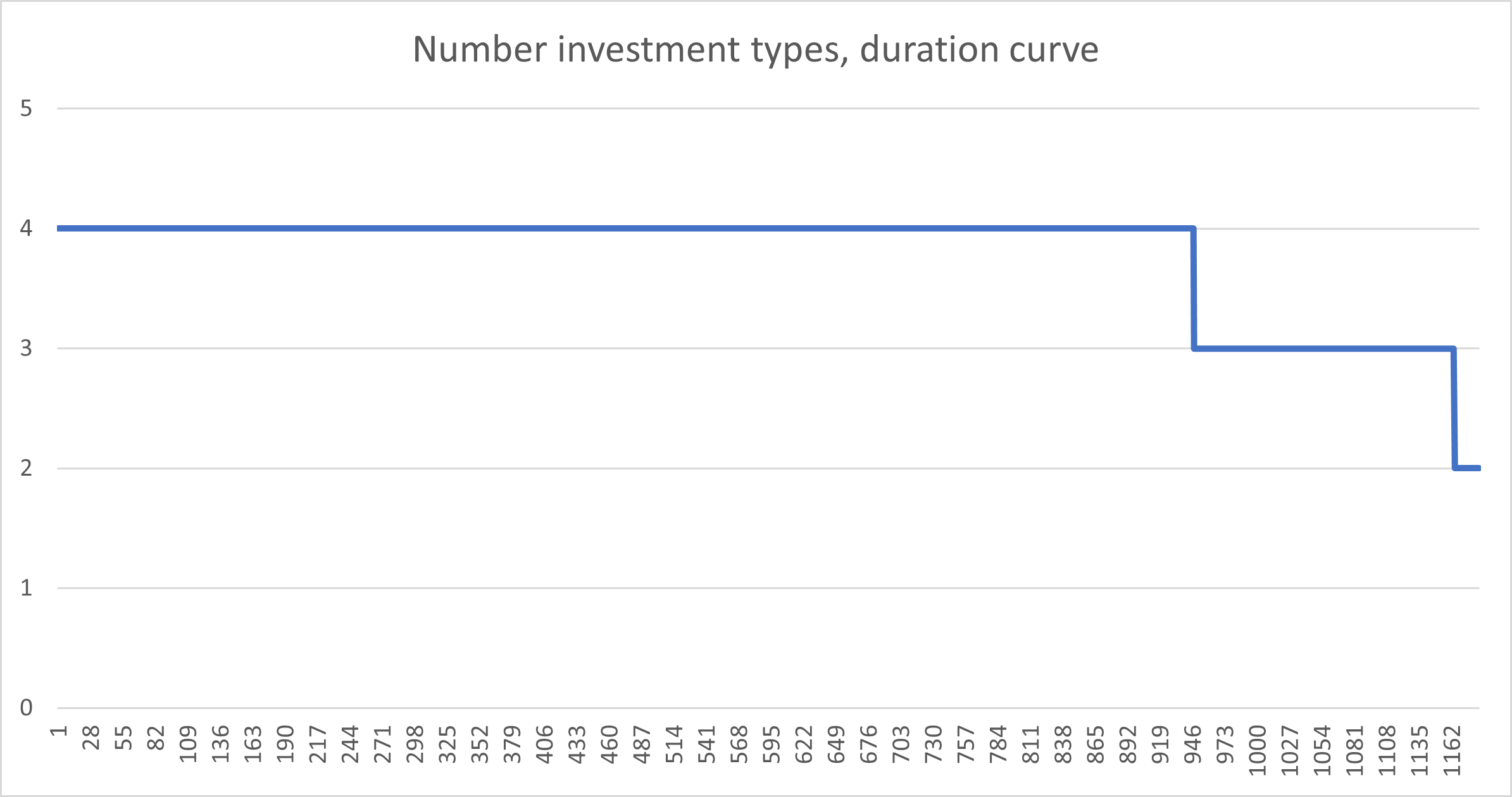}
    \caption{Duration curve of the number of investments types. Four types are considered: Skills, new technicians, overtime and digitization.}
    \label{fig:num_inv_types}
 \end{figure}

Finally, we note that investments are made in all four categories. A duration curve of the number of investments is illustrated in Figure~\ref{fig:num_inv} and a duration curve of how many of the four investment types are utilized is illustrated in Figure~\ref{fig:num_inv_types}. This highlight the method's capability in exploiting synergies between the investments instead of the more traditional approach where typically only one type of investment is analyzed at a time.

\begin{table}
\renewcommand{\tabcolsep}{0.75mm}
\scriptsize
    \centering
    \begin{tabular}{l|rrrrrrrrrrrr}\hline\hline
\myemph{Avg. number} &30 &45 &60 &75 &90  &105 &120 &135 &150 &210 &300 &600 \\
\myemph{new techs}   & tasks& tasks& tasks& tasks& tasks & tasks& tasks& tasks& tasks& tasks& tasks& tasks\\\hline
Narrow&1.10&1.74&2.76&3.10&3.20&3.20&3.20&3.20&3.76&3.00&8.40&18.24\\
Wide&0.28&0.62&1.40&2.14&2.78&2.96&3.08&3.18&3.49&3.00&7.88&15.53\\\hline
AllSkills&0.90&1.35&2.25&2.80&3.00&3.00&3.00&3.00&3.00&-&-&-\\
ConfigurationDeBase&0.85&1.35&2.25&2.80&3.00&3.00&3.00&3.00&3.92&3.00&8.14&16.96\\
nbrTech&0.30&0.70&1.85&2.70&3.45&3.85&4.00&4.00&4.00&-&-&-\\
RedSkills&0.90&1.35&2.25&2.80&3.00&3.00&3.00&3.00&3.00&-&-&-\\
TpsRep10-20&0.50&1.15&1.80&2.00&2.50&2.55&2.70&2.95&3.00&-&-&-\\
\hline\hline
\multicolumn{6}{l}{\textbf{Mathlouthi instances}}\klemme 
    \end{tabular}
    \caption{\label{tab:math_inv_tech}
    Avg. number of investments in new technicians}
\end{table}

\begin{table}
\renewcommand{\tabcolsep}{1.5mm}
\scriptsize
    \centering
    \begin{tabular}{l|rrrrr}\hline\hline
    \myemph{Avg. number new techs}&3 tech & 4 tech & 6 tech&12 tech&24 tech\\\hline
Narrow&2.70&3.03&5.96&12.00&23.24\\
Wide&2.17&2.49&5.48&11.35&19.05\\
\hline
AllSkills&2.48&-&-&-&-\\
ConfigurationDeBase&2.60&-&5.72&11.68&21.38\\
nbrTech&-&2.76&-&-&-\\
RedSkills&2.48&-&-&-&-\\
TpsRep10-20&2.13&-&-&-&-\\
\hline\hline
\multicolumn{6}{l}{\textbf{Mathlouthi instances}}\klemme 
    \end{tabular}
    \caption{\label{tab:math_inv_tech_techs}Avg. number of investments in new technicians}
\end{table}

\begin{table}
\renewcommand{\tabcolsep}{0.65mm}
\scriptsize
    \centering
    \begin{tabular}{l|rrrrrrrrrrrr}\hline\hline
\myemph{Avg. number}       &30 &45 &60 &75 &90  &105 &120 &135 &150 &210 &300 &600 \\
\myemph{skill investments} & tasks& tasks& tasks& tasks& tasks & tasks& tasks& tasks& tasks& tasks& tasks& tasks\\\hline
Narrow&5.08&7.22&8.84&11.10&12.44&13.22&15.16&15.60&17.06&15.30&39.22&81.76\\
Wide&4.06&5.04&7.00&8.76&12.00&12.60&13.56&15.16&16.24&16.00&38.08&77.78\\\hline
AllSkills&0.00&0.00&0.00&0.00&0.00&0.00&0.00&0.00&0.00&-&-&-\\
ConfigurationDeBase&5.30&7.55&9.45&12.20&14.25&14.65&15.75&16.50&17.88&15.65&38.65&79.88\\
nbrTech&5.65&8.05&10.80&12.85&16.10&17.70&20.65&22.75&22.30&-&-&-\\
RedSkills&7.30&9.40&12.70&15.70&19.15&20.55&22.00&23.70&23.45&-&-&-\\
TpsRep10-20&4.60&5.65&6.65&8.90&11.60&11.65&13.40&13.95&14.70&-&-&-\\
\hline\hline
\multicolumn{6}{l}{\textbf{Mathlouthi instances}}\klemme 
    \end{tabular}
    \caption{\label{tab:math_inv_skill}
    Avg. number of investments in skill upgrades.
Note that new technicians may also have their skills set upgraded and that each technician may have multiple upgrades.}
\end{table}

\begin{table}
\renewcommand{\tabcolsep}{1.5mm}
\scriptsize
    \centering
    \begin{tabular}{l|rrrrr}\hline\hline
    \myemph{Avg. number skill investments} &3 tech & 4 tech & 6 tech&12 tech&24 tech\\\hline
Narrow&10.83&16.43&26.32&58.90&95.52\\
Wide& 9.96&13.98&25.44&57.23&92.15\\\hline
AllSkills&0.00&-&-&-&-\\
ConfigurationDeBase&13.35&-&25.88&58.06&94.02\\
nbrTech&-&15.21&-&-&-\\
RedSkills&17.11&-&-&-&-\\
TpsRep10-20&10.12&-&-&-&-\\
\hline\hline
\multicolumn{6}{l}{\textbf{Mathlouthi instances}}\klemme 
    \end{tabular}
    \caption{\label{tab:math_inv_skill_techs}Avg. number of investments in skill upgrades.
Note that new technicians may also have their skills set upgraded and that each technician may have multiple upgrades.}
\end{table}

\begin{table}
\renewcommand{\tabcolsep}{0.6mm}
\scriptsize
    \centering
    \begin{tabular}{l|rrrrrrrrrrrr}\hline\hline
\myemph{Avg. number} &30 &45 &60 &75 &90  &105 &120 &135 &150 &210 &300 &600 \\
\myemph{digitized tasks} & tasks& tasks& tasks& tasks& tasks & tasks& tasks& tasks& tasks& tasks& tasks& tasks\\\hline
Narrow&1.52&2.24&3.98&5.64&8.52&11.04&14.26&19.70&24.40&63.70&47.32&84.69\\
Wide&4.70&6.32&8.72&11.06&14.04&17.12&20.48&25.52&28.56&61.50&57.74&109.33\\\hline
AllSkills&3.20&4.50&6.70&8.70&11.85&15.35&19.45&26.35&32.85&-&-&-\\
ConfigurationDeBase&3.20&4.50&6.70&8.70&11.85&15.45&19.30&26.40&26.56&62.60&52.53&96.28\\
nbrTech&3.20&4.50&6.25&8.05&10.80&12.70&15.30&18.80&22.80&-&-&-\\
RedSkills&3.20&4.50&6.70&8.85&11.75&15.50&19.45&26.40&32.70&-&-&-\\
TpsRep10-20&2.75&3.40&5.40&7.45&10.15&11.40&13.35&15.10&17.15&-&-&-\\
\hline\hline
\multicolumn{6}{l}{\textbf{Mathlouthi instances}}\klemme 
    \end{tabular}
    \caption{\label{tab:math_inv_dig}
    Avg. number of investments in task digitization}
\renewcommand{\tabcolsep}{0.4mm}
\end{table}

\begin{table}
\scriptsize
    \centering
    \begin{tabular}{l|rrrrr}\hline\hline
    \myemph{Avg. number digitized tasks} &3 tech & 4 tech & 6 tech&12 tech&24 tech\\\hline
Narrow&13.22&8.28&42.16&73.35&56.48\\
Wide&17.64&14.48&49.16&86.03&88.10\\
\hline
AllSkills&14.33&-&-&-&-\\
ConfigurationDeBase&21.47&-&45.66&79.69&70.53\\
nbrTech&-&11.38&-&-&-\\
RedSkills&14.37&-&-&-&-\\
TpsRep10-20&9.57&-&-&-&-\\
\hline\hline
\multicolumn{6}{l}{\textbf{Mathlouthi instances}}\klemme 
    \end{tabular}
    \caption{\label{tab:math_inv_dig_techs}Avg. number of investments in task digitization}
\end{table}

\begin{table}
\renewcommand{\tabcolsep}{0.7mm}
\scriptsize
    \centering
    \begin{tabular}{l|rrrrrrrrrrrr}\hline\hline
\myemph{Avg. number} &30 &45 &60 &75 &90  &105 &120 &135 &150 &210 &300 &600 \\
\myemph{overtime}    & tasks& tasks& tasks& tasks& tasks & tasks& tasks& tasks& tasks& tasks& tasks& tasks\\\hline
Narrow&2.70&2.84&3.72&4.52&4.78&5.36&5.48&5.86&6.57&6.00&15.02&30.20\\
Wide&1.78&2.58&3.30&4.50&5.04&5.38&5.42&5.58&6.50&5.90&15.38&30.63\\\hline
AllSkills&2.55&3.15&3.70&4.80&5.15&5.65&5.50&5.80&5.70&-&-&-\\
ConfigurationDeBase&2.55&3.15&3.80&4.85&5.15&5.55&5.50&5.80&7.12&5.95&15.20&30.40\\
nbrTech&2.40&2.85&3.80&5.00&5.65&6.70&7.05&7.25&7.55&-&-&-\\
RedSkills&2.55&3.15&3.70&4.90&5.10&5.50&5.50&5.80&5.70&-&-&-\\
TpsRep10-20&1.15&1.25&2.55&3.00&3.50&3.45&3.70&3.95&4.25&-&-&-\\
\hline\hline
\multicolumn{6}{l}{\textbf{Mathlouthi instances}}\klemme 
    \end{tabular}
    \caption{\label{tab:math_inv_over}
    Avg. number of investments in overtime.
Note that new technicians may also be assigned overtime.}
\end{table}

\begin{table}
\renewcommand{\tabcolsep}{1.5mm}
\scriptsize
    \centering
    \begin{tabular}{l|rrrrr}\hline\hline
    \myemph{Avg. number overtime}    &3 tech & 4 tech & 6 tech&12 tech&24 tech\\\hline
Narrow&4.45&5.46&10.92&21.83&35.16\\
Wide&4.27&5.27&10.70&22.35&37.45\\\hline
AllSkills&4.67&-&-&-&-\\
ConfigurationDeBase&4.94&-&10.81&22.09&36.18\\
nbrTech&-&5.36&-&-&-\\
RedSkills&4.66&-&-&-&-\\
TpsRep10-20&2.98&-&-&-&-\\
\hline\hline
\multicolumn{6}{l}{\textbf{Mathlouthi instances}}\klemme 
    \end{tabular}
    \caption{\label{tab:math_inv_over_techs}Avg. number of investments in overtime.
Note that new technicians may also be assigned overtime.}
\end{table}

\bigskip

\begin{figure}
\centering
    \includegraphics[width=0.8\textwidth]{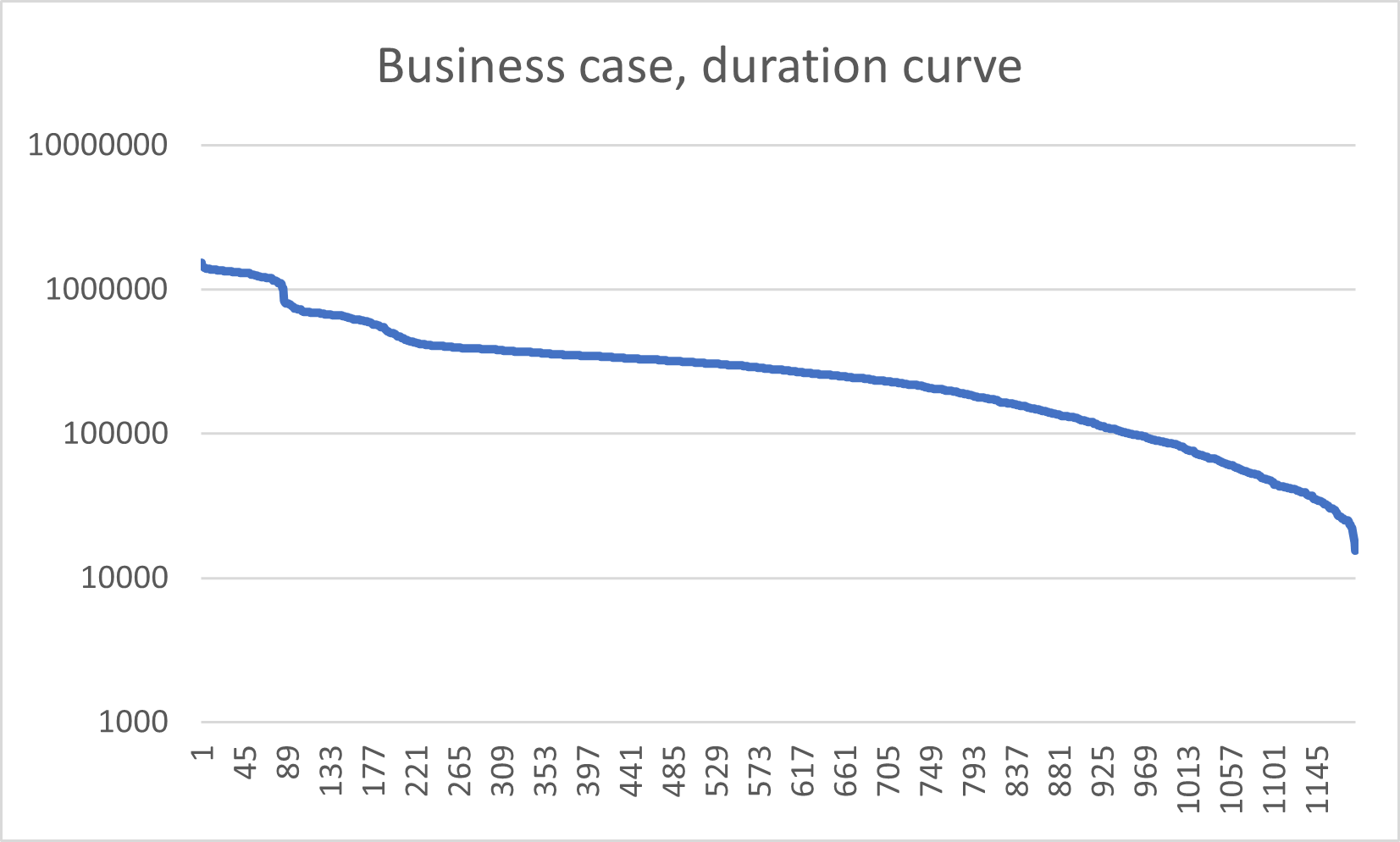}
    \caption{Duration curve of the business case value for the investments for the Mathlouthi instances. The business case is calculated as the savings in the ALNS with and without investments, subtracted the investment costs: ALNS(inv) - ALNS(no inv) - CAPEX. Note the logarithmic y axis.}
    \label{fig:bc_log_scale}
 \end{figure}

The investments are added to the original Mathlouthi instances and then solved by the ALNS. The purpose  is to analyze the consequence of the investments in the full TRSP. The first comparison concerns the actual business case; does the OPEX saving exceed the CAPEX cost? The OPEX saving is calculated as the ALNS OPEX without investments minus the ALNS OPEX with investments. The ALNS OPEX is the objective function in the TRSP \eqref{trsp_obj}. The CAPEX is the total cost of the investments, decided by the ASSM. Figure~\ref{fig:bc_log_scale} shows a duration curve for the business case values and the average business case values are found in tables \ref{tab:math_inv_bc} and \ref{tab:math_inv_bc_techs}. The business case for each instance is positive and ranges from appr. 15000 to 1500000. Recall that the objective function value is given in travel time, see Section~\ref{investment_data}. Reverting to actual cost, the business case range is from appr. \xeuro 3000 to appr. \xeuro 300000. 

\begin{figure}
\centering
    \includegraphics[width=0.8\textwidth]{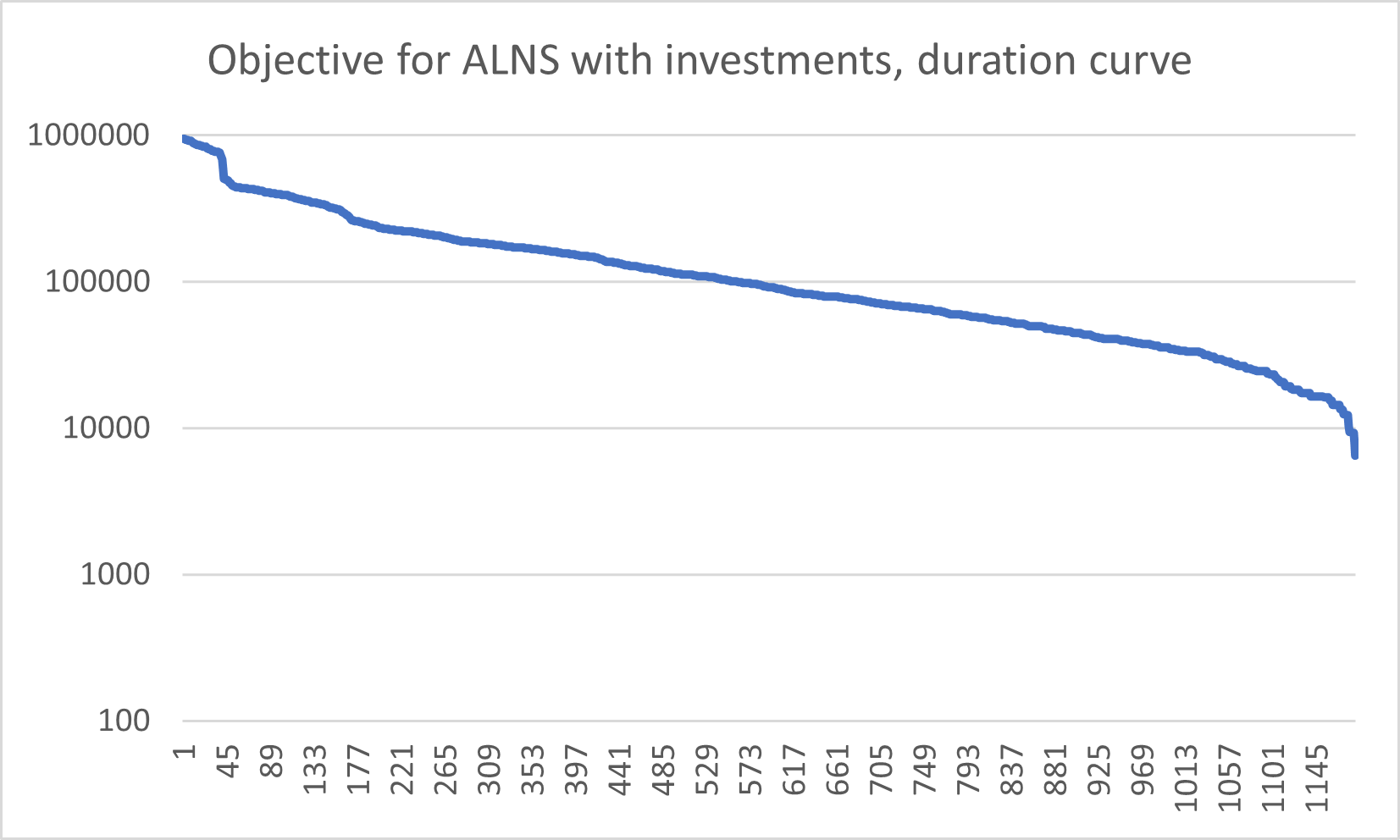}
    \caption{Duration curve of the OPEX of the ALNS for the Mathlouthi instances including investment decision. Note the logarithmic y axis.}
    \label{fig:opex_log_scale}
 \end{figure}

A duration curve of the OPEX for the Mathlouthi instances with investments are shown in Figure~\ref{fig:opex_log_scale}. They span from appr. 6500 to appr. 950000 which corresponds to appr. \xeuro 1300 to \xeuro 190000. The average saving in OPEX is displayed in tables \ref{tab:math_inv_opex_diff} and \ref{tab:math_inv_opex_diff_techs} spanning from 40000 to appr. 1575000 which is \xeuro 8000 to \xeuro 3150000. The saving in OPEX is significant.

\begin{figure}
\centering
    \includegraphics[width=0.8\textwidth]{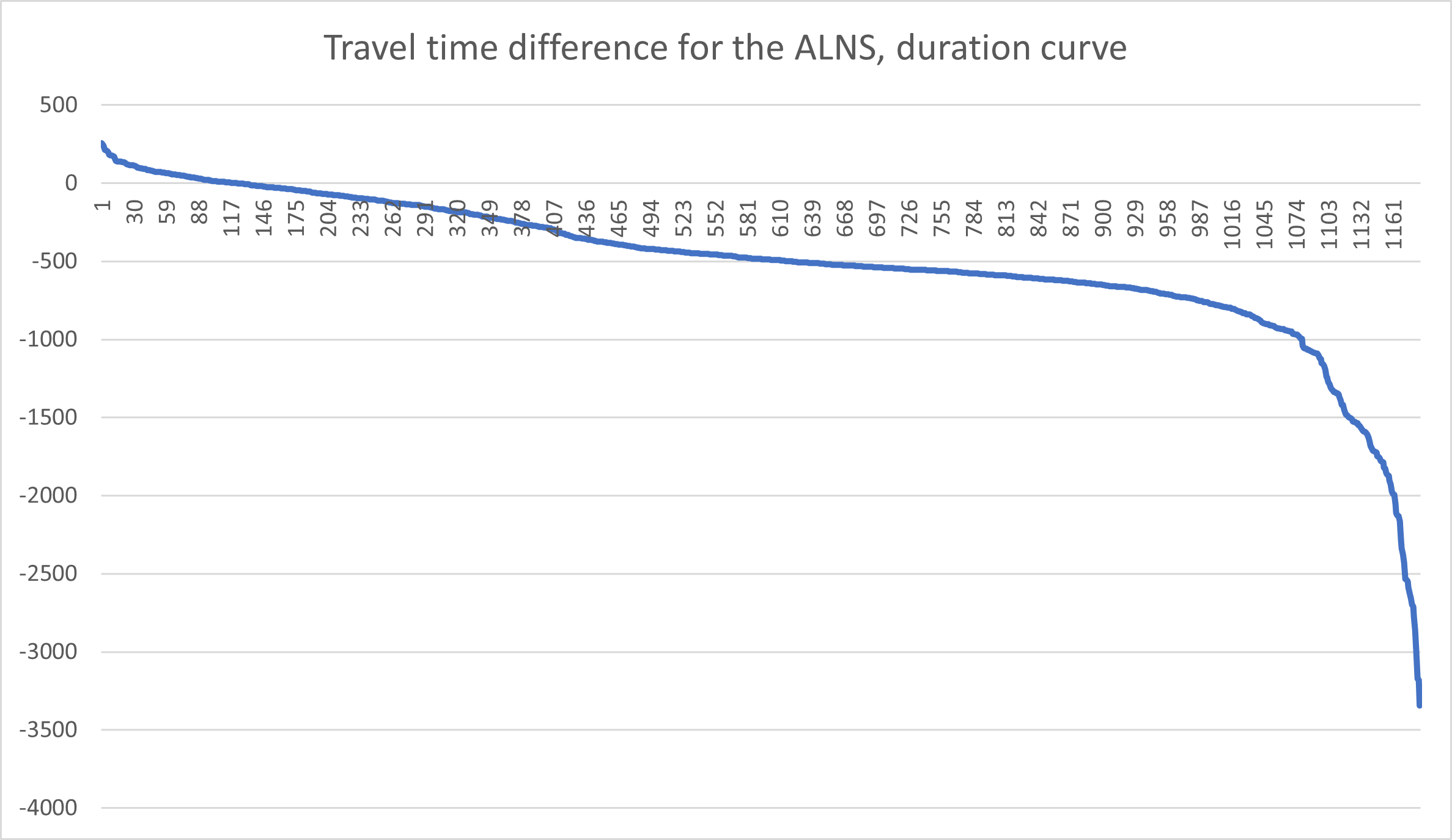}
    \caption{Duration curve of change in the travel time for the ALNS with investments and the ALNS without investments, calculated as ALNS(no inv) - ALNS(inv).}
    \label{fig:alns_travel_time_diff}
 \end{figure}

The next investigation is the driving force behind the investments. We compare the results for the ALNS with investments and without investments. The average difference in the number of unserved tasks is provided in tables \ref{tab:math_inv_unserved_tasks} and \ref{tab:math_inv_unserved_tasks_techs}, and a duration curve of the difference in travel time is shown in Figure~\ref{fig:alns_travel_time_diff}. The results clearly show that reducing the number of unserved tasks drives investments. In fact, the travel time tends to increase after investments are made. The reason is two-fold; the penalty of unserved tasks dominate the objective function and is given higher priority than travelling. Also, serving more tasks very likely will increase the travel time because more locations are visited.

The OPEX savings suggest that the business cases are generally very strong. Of course, the instances only span a single day which is not enough to make actual investment decisions. Still, the results strongly suggest that the OPEX for the Mathlouthi instances could be optimized significantly through strategic decisions.

\begin{table}
\renewcommand{\tabcolsep}{1.0mm}
\scriptsize
    \centering
    \begin{tabular}{l|rrrrrrrrrrrr}\hline\hline
\myemph{Avg. business case} &30 tasks&45 tasks&60 tasks&75 tasks&90 tasks &105 tasks\\\hline
Narrow&39572.56&71932.06&121807.27&168483.44&217591.93&263747.35\\
Wide&43036.02&65584.00&106535.19&156765.33&213806.93&263026.00\\\hline
AllSkills&36997.16&61877.78&105740.98&160539.50&212046.82&257545.38\\
nbrTech&41845.71&68803.14&109106.42&157494.90&212332.63&266277.46\\
RedSkills&48744.99&85586.63&136325.97&187695.97&236751.35&280603.18\\
TpsRep10-20&31262.52&46402.15&88637.97&125027.61&183053.74&233725.72\\
ConfigurationDeBase&47671.09&81120.43&131044.82&182363.97&234312.62&278781.64\\ \hline \hline
&120 tasks&135 tasks&150 tasks&210 tasks&300 tasks&600 tasks\\\hline
Narrow&300898.46&348778.73&371975.63&506711.24&667242.94&1314148.90\\
Wide&306920.51&355101.67&370678.97&501449.06&665625.99&1259082.05\\\hline
AllSkills&291524.26&337652.90&368729.79&-&-&-\\
nbrTech&310405.40&353992.98&357348.18&504080.15&666434.47&1288235.09\\
RedSkills&313753.82&362059.74&394572.71&-&-&-\\
TpsRep10-20&314995.41&354725.06&385868.29&-&-&-\\
ConfigurationDeBase&288868.53&351270.32&406034.00&-&-&-\\
\hline\hline
\multicolumn{6}{l}{\textbf{Mathlouthi instances}}\klemme 
    \end{tabular}
    \caption{\label{tab:math_inv_bc}
    The average business case value}
\end{table}

\begin{table}
\scriptsize
\renewcommand{\tabcolsep}{1.5mm}
    \centering
    \begin{tabular}{l|rrrrr}\hline\hline
    \myemph{Avg. business case}&3 tech & 4 tech & 6 tech&12 tech&24 tech\\\hline
Narrow&227146.87&215208.96&548953.51&1000010.01&1291207.76\\
Wide&226447.56&212845.82&548819.16&973789.48&1184787.57\\\hline
AllSkills&203628.28&-&-&-&-\\
ConfigurationDeBase&268904.79&-&548886.33&986899.75&1243909.90\\
nbrTech&-&214027.39&-&-&-\\
RedSkills&225699.65&-&-&-&-\\
TpsRep10-20&194920.28&-&-&-&-\\
\hline\hline
\multicolumn{6}{l}{\textbf{Mathlouthi instances}}\klemme 
    \end{tabular}
    \caption{\label{tab:math_inv_bc_techs}The average business case value}
\end{table}

\begin{table}
\renewcommand{\tabcolsep}{1.3mm}
\scriptsize
    \centering
    \begin{tabular}{l|rrrrrrrrrrrr}\hline\hline
\myemph{Avg. business case} &30 tasks&45 tasks&60 tasks&75 tasks&90 tasks &105 tasks\\\hline
Narrow&46085.36&81150.76&137052.67&188725.94&245318.33&298062.05\\
Wide&56065.12&83465.40&131745.19&186985.33&254930.93&312240.00\\\hline
AllSkills&47224.66&76165.28&126855.98&186392.00&247589.32&302062.88\\
ConfigurationDeBase&58024.09&95672.18&152535.57&208648.47&270353.87&324016.89\\
nbrTech&51483.46&82457.39&129039.42&181770.40&246578.63&306281.96\\
RedSkills&59227.99&100203.13&157885.47&215935.47&272691.60&326147.43\\
TpsRep10-20&39416.02&57042.40&105678.22&146531.86&213409.74&267245.97\\\hline\hline
&120 tasks&135 tasks&150 tasks&210 tasks&300 tasks&600 tasks\\\hline
Narrow&343385.06&405051.73&441034.24&672796.74&803754.64&1564215.90\\
Wide&363497.01&425759.27&449748.08&662014.06&827685.79&1567527.92\\\hline
AllSkills&346224.26&409737.90&457019.79&-&-&- \\
ConfigurationDeBase&365281.65&426780.48&432281.98&667405.40&815720.22&1565774.50\\
nbrTech&360699.07&417918.49&460550.71&-&-&-\\
RedSkills&367382.66&427764.56&474604.04&-&-&-\\
TpsRep10-20&327617.53&394826.07&454936.00&-&-&-\\
\hline\hline
\multicolumn{6}{l}{\textbf{Mathlouthi instances}. Gap $=$ ALNS(no inv) - ALNS(inv).}\klemme 
    \end{tabular}
    \caption{\label{tab:math_inv_opex_diff}
    The average difference between the OPEX for the ALNS without and with investments}
\end{table}

\begin{table}
\renewcommand{\tabcolsep}{1.5mm}
\scriptsize
    \centering
    \begin{tabular}{l|rrrrr}\hline\hline
    \myemph{Avg. business case}&3 tech & 4 tech & 6 tech&12 tech&24 tech\\\hline
Narrow&265815.90&242573.57&667340.71&1209667.76&1479460.96\\
Wide&275057.67&254488.55&684000.56&1214532.35&1447975.32\\\hline
AllSkills&244363.56&-&-&-&-\\
ConfigurationDeBase&328268.02&-&675670.63&1212100.06&1465467.34\\
nbrTech&-&248531.06&-&-&-\\
RedSkills&266871.37&-&-&-&-\\
TpsRep10-20&222967.09&-&-&-&-\\
\hline\hline
\multicolumn{6}{l}{\textbf{Mathlouthi instances}. Gap $=$ ALNS(no inv) - ALNS(inv).}\klemme 
    \end{tabular}
    \caption{\label{tab:math_inv_opex_diff_techs}The average difference between the OPEX for the ALNS without and with investments}
\end{table}

\begin{table}
\renewcommand{\tabcolsep}{0.5mm}
\scriptsize
    \centering
    \begin{tabular}{l|rrrrrrrrrrrr}\hline\hline
               &30 &45 &60 &75 &90  &105 &120 &135 &150 &210 &300 &600 \\
\myemph{Unserved tasks} & tasks& tasks& tasks& tasks& tasks & tasks& tasks& tasks& tasks& tasks& tasks& tasks\\\hline
Narrow&8.32&14.84&26.48&35.24&43.26&49.64&55.90&63.98&71.06&101.10&143.62&292.87\\
Wide&9.42&14.66&25.08&36.12&48.18&55.24&61.62&69.48&75.02&99.90&149.54&291.45\\\hline
AllSkills&7.80&13.75&25.65&36.55&45.25&50.30&55.05&63.40&69.30&-&-&-\\
ConfigurationDeBase&10.50&17.55&29.35&39.10&47.40&52.40&57.00&64.35&72.73&100.50&146.58&292.20\\
nbrTech&8.60&13.95&24.05&34.80&47.20&56.05&63.70&68.35&73.00&-&-&-\\
RedSkills&10.70&18.05&29.75&39.80&47.35&52.75&57.00&65.30&70.90&-&-&-\\
TpsRep10-20&6.75&10.45&20.10&28.15&41.40&50.70&61.05&72.25&80.50&-&-&-\\
\hline\hline
\multicolumn{10}{l}{\textbf{Mathlouthi instances.} Gap $=$ ALNS(no inv) - ALNS(inv).}\klemme 
    \end{tabular}
    \caption{\label{tab:math_inv_unserved_tasks}
    The average number of unserved tasks by the ALNS without investments and the ALNS with investments}
\end{table}

\begin{table}
\renewcommand{\tabcolsep}{1.5mm}
\scriptsize
    \centering
    \begin{tabular}{l|rrrrr}\hline\hline
    \myemph{Unserved tasks}&3 tech & 4 tech & 6 tech&12 tech&24 tech\\\hline
Narrow&43.70&41.37&112.16&213.25&309.76\\
Wide&46.60&45.23&118.50&219.10&297.15\\\hline
AllSkills&40.78&-&-&-&-\\
ConfigurationDeBase&52.57&-&115.33&216.18&304.16\\
nbrTech&-&43.30&-&-&-\\
RedSkills&43.51&-&-&-&-\\
TpsRep10-20&41.26&-&-&-&-\\
\hline\hline
\multicolumn{6}{l}{\textbf{Mathlouthi instances.} Gap $=$ ALNS(no inv) - ALNS(inv).}\klemme 
    \end{tabular}
    \caption{\label{tab:math_inv_unserved_tasks_techs}The average number of unserved tasks by the ALNS without investments and the ALNS with investments}
\end{table}

\subsubsection{TDC Net instance}
Preliminary test runs reveal that adding investment decisions to the TDC Net instance makes the problem significantly harder to solve. We thus warmstart the ASSM with the solution without investment options from Section~\ref{results_tdc_bc}. Other than that, the ASSM and ALNS are run with the same parameter settings as without investments in Section~\ref{results_tdc_bc}.

Results for the ASSM and the ALNS are summarized in Table~\ref{tab:tdc_investment_results}. They show that the ASSM times out after 16 column generation iterations and that solving the final integer master problem also times out. It generates 8545 paths in total. Recall that without investments, the ASSM managed 38 iterations and generated 19484 before time out. This indicates that including investment decisions makes the problem harder to solve. The reason is twofold; the subproblem becomes more time consuming because of the investment variables; and the relaxed master problem also becomes more time consuming to solve also because of the investment variables. The ALNS also manages fewer iterations and generated paths, probably because the number of technicians has increased. 

We observe a difference in the number of unserved paths and also a difference in the objective function value, just like we did for the TDC Net instance without investments. If we subtract the penalty of unserved tasks from the objective, we again see that the travel times are quite different, with the ASSM having a large estimated travel time.

The actual investments found by the ASSM are summarized in Table~\ref{tab:tdc_investments}. The number of investments seems fair compared to the instance size. With  6 new technicians, the total number of daily technicians in the instance is 745 (5 days of 6 new technicians). 95 skills upgrade take place and 35 technicians must work overtime. Also, of the 2677 tasks, 7 are digitized. 

The total investment costs sum up to 43075, which is a small value compared to the objective function value. Using the ALNS results, we can quantify the business case for the TRSP as 441858 - 183977 - 43075 = 214806. The business case is clearly positive, which indicates attractive investment suggestions. The TDC Net instance only spans 5 days, so more analyses would be needed to ensure the robustness of the investments. 

A closer look at the ALNS solutions without and with investments, show that the largest decrease in the objective stems going from fewer unserved tasks, from 42 to 16 unserved tasks. It can be discussed whether a solving 26 more tasks justify employment of 6 new technicians.
TDC Net is very focused on customer satisfaction and thus prioritize task service very highly. Still, if the number of new technicians is too high, then the penalty of unserved tasks should be lowered. The higher penalty, the more investments. It is, however, interesting to see that the amount of overtime does not explode but instead the investments indicate that new technicians and skill upgrades are an attractive alternative to increase customer satisfaction.

Unlike for the Mathlouthi instances, we do see a significant decrease in the travel time in the ALNS when including investments. Without investments the travel time was 35501 (objective value minus penalty: 441858-406357) and with investments 29813 (183977-154164). This is a reduction of 5688 travel minutes or 16\%! The saving in travel time corresponds to almost 95 hours of driving in just 5 days. Relative to the total number of tasks, the reduction in the number of unserved tasks is very small and thus does not have a negative impact on the travel time as was the case for the Mathlouthi instances.

Still, the decrease in travel time is smaller than the investment costs. The investments are driven by unserved tasks as the penalty dominates the objective function values. Had the focus been on investments to reduce travel time, more emphasis should be put on improving the travel time estimate in the ASSM. Also, if the ASSM scaled better, we would probably be able to identify even better investments. In the following section, we discuss future research directions concerning these aspects.

\begin{table}[htbp]
\scriptsize
    \centering
    \begin{tabular}{l|rr}\hline\hline
    \myemph{Comparison}
        & ASSM & ALNS\\\hline
        Objective & 284678&183977\\
        Total time sec & 11341.65&14039.03\\
        \# unserved tasks & 6& 16\\
        Total penalty & 47668 &154164\\
        Avg. path length & 6.14 &4.67\\
        \# paths in solution & 560 &568 \\
        Total number generated paths & 8545 & 271522\\\hline
        Number cg. iters. & 16 & - \\
        Optimal cg. & False & - \\
        Optimal int. master problem & False & - \\ \hline
        Avg. number improvements & - &  7\\
        Avg. number iters. & - &  1200 \\
        Avg. last impr. iter. & - & 932 \\ \hline \hline
\multicolumn{3}{l}{\textbf{TDC Net instances}}\klemme 
  \end{tabular}
    \caption{Results for solving the 
instance with the ASSM and the ALNS with investments. The last two sections concern results for the ASSM resp. the ALNS only.}
    \label{tab:tdc_investment_results}
\end{table}

\begin{table}[htbp]
\scriptsize
    \centering
    \begin{tabular}{l|r}\hline\hline
    \myemph{Investments} & \\ \hline
    \# skill upgrades & 95 \\
    Skill upgr. cost & 16625 \\
    \# new techs & 6 \\
    \# new daily techs & 30 \\
    New tech cost & 7200\\
    \# dig. tasks & 7 \\
    Dig. task cost & 3500 \\
    \# overtime & 35 \\
    Overtime cost & 15750\\ \hline \hline
\multicolumn{2}{l}{\textbf{TDC Net instances}}\klemme 
  \end{tabular}
    \caption{The actual investments found by the ASSM}
    \label{tab:tdc_investments}
\end{table}

\section{Future work}
\label{future}
The investment master problem is very generic and could be applied to many other problems than the TRSP. Most obvious is to solve other routing problems, but generally any set cover like problem could be considered.
In this section, however, we remain focused on the TRSP and possible extensions of the proposed solution approach.

\subsection{Investments in the TRSP}
If investments are motivated by reducing the travel time and not the number of unserved tasks, then it would be beneficial to consider investments in the original TRSP. That is, instead of approximating the TRSP with the task assignment approach, the full TRSP should be considered. This would eliminate the inaccuracies from estimating the travel time in the assignment approach. The column generation approach could be preserved with the same master problem and where the subproblem becomes the TRSP for a single technician. To limit the size of of the subproblem, an upper bound can be set on how far technicians can travel from their home depot. Still, the subproblem remains difficult to solve and it is fair to assume that the approach would require significantly more time to find a solution of acceptable quality.

\subsection{Reducing the problem instance size}
To obtain robust investments, it is desirable to solve the investment problem for a larger time horizon, e.g., a year. This, however, would lead to very large problem instances with many daily technicians (and subproblems) and it is fair to assume the column generation approach would require many iterations before the investments converge.

The time horizon can be reduced by clustering similar time periods, as done for the Capacity Expansion Problem (CEP), see Section~\ref{literature} and the work of Buchholz et al. \cite{buchholz}. Specifically, k-means clustering can be applied to cluster weeks according to task types, number of technicians and the total task duration. We suggest to cluster weeks instead of days, as many tasks have time windows spanning more than a day. If this is not the case in certain instances, then clustering days may be a better approach as it is probably possible to reduce the time horizon even more.
%

\bigskip

Another method to reduce the problem instance size is to reduce the geographical scope of the problem. Obviously, this is a heuristic approach, especially if the instance contains tasks or technicians close to the border between geographical areas.

Decomposing problems according to geography is a known method from the VRP research area. Santin et al. \cite{vigo} propose and survey decomposition methods for heuristics solving the VRP. They argue that decomposing the VRP into independent subproblems is beneficial as this allows for parallelism. Decompositions in the literature are based on geography \cite{Groer}; on vehicle capacity \cite{vigo}; randomly \cite{vigo}; using unsupervised machine learning such as k-means clustering \cite{arthur} or historical related k-medoids clustering \cite{park}; or by collapsing tasks into hyper nodes in a path based decomposition \cite{vigo}. Santin et al. \cite{vigo} evaluate the decomposition methods on instances with 600 customers and the unsupervised machine learning methods show best results. The TRSP is a rich VRP, so it is fair to assume that results from the VRP would also give attractive results for the TRSP.

\subsection{Investment in new technicians}
The location of new technicians can be optimized by setting the same travel time from their depot to every task. The estimated travel time of a resulting path would then be set by the travel time between two assigned tasks. This approach would, however, make the subproblem more difficult to solve, because all tasks are equal candidates to be assigned to the technician. This could be remedied by clustering tasks according to their geographical location and generate a new technician for each cluster with equal travel time to the tasks in that cluster.

\subsection{Limiting the investments}

For practical reasons it may not be possible to make too many changes to the technician staff. For instance,
we cannot upgrade skills of all technicians at the same time due to capacity limits. Such constraints can easily be added to the model as budget constraints for each investment, or an overall investment pool budget.


\section{Conclusion}
\label{conclusion}
This paper proposed a method for building scenarios for the Technician Routing and Scheduling Problem (TRSP). Provided is a base case with a certain set of tasks to serve, technician fleet, skills set and work schedule. The proposed method is capable of building an investment scenario such that the reduction in TRSP costs exceeds the investment costs. The investigated investments were; employing new technicians, upgrading the skill sets of each technician, assigning overtime to each technician and digitizing the task equipment allowing for remote repairs. Each of these investments have an associated cost. 

In real-life, scenarios and strategic decisions are often based on what-if analyses where single investments are investigated one at a time. The proposed method excels in its holistic view and capability of investigating synergies between several investments. To the best of our knowledge, such an approach has not been proposed previously to the TRSP and similar approaches for the related Vehicle Routing Problem have a much more narrow scope to e.g. optimizing the fleet. We thus believe that the proposed method helps close a research gap on strategic decisions in a routing context.

The proposed method consists of a matheuristic based on column generation. The master problem minimizes the total travel time, penalties of unserved tasks and investment costs. The subproblem assigns tasks to each technician subject to possible investments. Instead of calculating the full routing problem for the technician, it instead estimates the travel time. This reduces the solution time of the subproblem. 

The paper furthermore proposed an ALNS meta-heuristic for the TRSP to evaluate if the investments found by the matheuristic are also beneficial to the full TRSP.
The ALNS was compared to a MIP model for the TRSP on smaller instances and results were very satisfactory with the ALNS finding solutions close to the MIP solutions.

The proposed matheuristic was evaluated on benchmark instances from the literature and on a real-life test instance from the telecommunication company TDC Net. Every run of the matheuristic resulted in scenarios that also produced a positive business case in the ALNS. That is, the reduced TRSP costs exceeded the investment costs. A closer look at the results revealed that the investments were driven by reducing the number of unserved tasks. The penalty of an unserved task is typically very high and thus dominates the objective function. In the real-life test instance, travel time was also decreased with around 16\%. 

The proposed method of estimating travel time is reasonable in this context but if the driving force of investments is to reduce travel time, then it would be more appropriate to adapt the matheuristic to solve the full TRSP. This would, however, be more time consuming.

The paper concluded with proposing future research directions to further the applicability of the matheuristic. The research directions both included suggestions on how to support other relevant investments and scenarios and on how to reduce the runtime of the method.

\section*{Acknowledgements}
We thank TDC Net for providing data and useful information on the real-life technician routing and scheduling problem. 
We also thank Stefan Røpke for sharing his ALNS code.

This work was funded by the Innovation Fund Denmark as part of the Grand Solution research project, GREENFORCE, project number 0224-00055B.

\bibliographystyle{plain}
\bibliography{bibliography}

\end{document}